\documentclass[leqno,final]{siamltex}
\usepackage{amsmath}
\usepackage{graphicx}
\usepackage{amssymb}

\newcommand{\hP}{\widehat{P}}

\newcommand{\DG}{\mbox{\tiny DG}}
\newcommand{\cE}{\mathcal{E}}
\newcommand{\cT}{\mathcal{T}}

\newcommand{\mce}{\mathcal{E}_h}
\newcommand{\mct}{\mathcal{T}_h}
\newcommand{\bl}{\big\langle}
\newcommand{\br}{\big\rangle}

\newcommand{\eps}{\epsilon}

\newcommand{\Ome}{\Omega}
\newcommand{\p}{\partial}
\newcommand{\nab}{\nabla}

\def\esssupI{\underset{t\in [0,\infty)}{\mbox{\rm ess sup }}}
\begin{document}

\title{Analysis of a second order discontinuous Galerkin  finite element method for the Allen-Cahn equation and the curvature-driven geometric flow}
\markboth{HUANRONG LI AND JUNZHAO HU}{2ND ORDER DG METHODS FOR ALLEN-CAHN EQUATION}

\author{
Huanrong Li\thanks{College of Mathematics and Statistics,
Chongqing Technology and Business University, Chongqing 400067,
China. ({\tt lihuanrong1979@163.com.})
The work of this author was  partially supported by National Science Foundation of China(11101453),
 Natural Science Foundation Project of Chongqing  CSTC(2013jcyjA20015, 2015jcyjA00009), and Chongqing Education Board of Science Foundation( KJ1400602) . Corresponding Author.}
\and
Junzhao Hu\thanks{Department of Mathematics, Iowa State University,
Ames, IA 50011, U.S.A. ({\tt junhu@iastate.edu.})
}
}

\maketitle

\begin{abstract}
The paper proposes and analyzes an efficient second-order in time numerical approximation for the Allen-Cahn equation, which is a second order nonlinear equation arising from the phase separation model. We firstly present a fully discrete interior penalty discontinuous Galerkin (IPDG) finite element method, which is based on the modified Crank-Nicolson scheme and a mid-point approximation of the nonliner term $f(u)$.  We then derive the  stability analysis and  error estimates for the proposed IPDG method under some regularity assumptions on the initial function $u_0$.  There are two key works in our analysis, one is to establish unconditionally energy-stable scheme for the discrete solutions. The other is to use a discrete spectrum estimate to handle the midpoint of the discrete solutions $u^m$ and $u^{m+1}$ in the nonlinear term, instead of using the standard Gronwall inequality technique. This discrete spectrum estimate is not trivial to obtain since the IPDG space and the conforming $H^1$ space are not contained in each other. We obtain that all our error bounds depend on  reciprocal of the perturbation parameter $\epsilon$  only in some lower polynomial order, instead of exponential order. These sharper error bounds are the key elements in proving the convergence of our numerical solution to the mean curvature flow. Finally, numerical experiments are also provided to show the performance of  the presented approach and method.
 \end{abstract}

\section{Introduction}\label{sec-1}
Let $\Omega\subseteq R^{d}(d=2,3)$ be a bounded polygonal or polyhedral domain. Consider the following nonlinear singular perturbation  model of the reaction-diffusion equation
\begin{equation}\label{eq1.1}
u_t-\Delta u+\frac{1}{\epsilon^2}f(u)=0,  \qquad \mbox{in }
\Omega_T:=\Omega\times(0,T).
\end{equation}
And we consider the following homogenous
Neumann boundary condition
\begin{equation}\label{eq1.2}
\frac{\partial u}{\partial \mathbf{n}} =0, \qquad \mbox{in }
\partial\Omega_T:=\partial\Omega\times(0,T),
\end{equation}
and initial condition
\begin{equation}\label{eq1.3}
u =u_0, \qquad \mbox{in }\Omega\times\{t=0\},
\end{equation}
where, $\mathbf{n}$ denotes the unit
outward normal vector to the boundary $\partial\Omega$, and the boundary condition \eqref{eq1.3} means that no mass loss occurs through the boundary walls.

Equation \eqref{eq1.1}, which is called  the Allen-Cahn equation, was originally introduced by Allen and Cahn in \cite{1, 8, 12} to  describe an interface evolving in time in the phase separation process of the crystalline solids. Herein, $\epsilon>0$ is a parameter related to the interface thickness, which is small compared to the characteristic length of the laboratory scale. $u$ denotes the concentration
of one of the two metallic species of the alloy, and $f(u)=F'(u)$  with
$F(u)$ being some given energy potential. Several choices of $F(u)$ have been presented in the literature \cite{2,9,3,4,5,6}. In this paper we focus on the following Ginzburg-Landau double-well potential \cite{feng2014finite, feng2017finite}
\begin{equation}\label{eq1.4}
F(u)=\frac{1}{4}(u^2-1)^2\
\  \mathrm{and}\  \ f(u)=F'(u)=(u^{2}-1)u.
\end{equation}
Although the potential term \eqref{eq1.4}  has been widely used, its quartic growth at
infinity leads to a variety of technical difficulties in the numerical approximation for the  Allen-Cahn equation. For example, in order to assure that our numerical scheme is second-order in time, we have to employ the modified Crank-Nicolson scheme and a  second order in time approximation of the potential term $f(u)$(see (3.4) in section 3.1).

An important feature of the Allen-Cahn equation (1.1) is that it can be viewed as the gradient flow with the Liapunov energy functional
\begin{equation}\label{eq1.5}
\textit{J}_\epsilon(u):=\int_\Omega \phi_\epsilon(u)dx \qquad \mbox{and}\qquad \phi_\epsilon(u)=\frac{1}{2}|\nabla u|^2+\frac{1}{\epsilon^2}F(u).
\end{equation}
More precisely, by taking the inner product of (1.1) with $-\Delta u+\frac{1}{\epsilon^2}f(u)$,
we immediately get the following energy law for (1.1)
\begin{equation}\label{eq1.6}
\frac{\partial}{\partial t}\textit{J}_\epsilon(u(t))=-\int_\Omega|-\Delta u+\frac{1}{\epsilon^2}f(u)|^2dx.
\end{equation}

Nowadays, the Allen-Cahn equation has  been  extensively investigated
 due to its connection to the interesting and complicated \emph{curvature-driven geometric flow} known as
\emph{the mean curvature flow} or\emph{ the motion by mean curvature} (cf.\cite{25, 8} and the references therein). It was proved that(see \cite{25}),
as $\epsilon\rightarrow 0$, the zero
level set of the solution $u$ of the problem (1.1)-(1.4), denoted by $\Gamma_t^\epsilon:=\{x\in\Omega; u(x,t)=0\}$ converges to the curvature-driven geometric flow as $\epsilon$ and mesh sizes $h$ and $k$ all tend to
zero, which refers to the evolution of a surface governed by the geometric law $V=\kappa$,
where $V$ is the inward normal velocity of the surface $\Gamma_t$ and $\kappa$ is its mean curvature,
see \cite{1,9}.

The Allen-Cahn  equation has been widely used in many complicated moving interface problems in fluid dynamics, materials science, image processing and biology (cf.\cite{12,Feng_Li15} and the references therein). Therefore, it is very important to develop
accurate and efficient numerical schemes to solve the Allen-Cahn equation. There are several challenges to obtain numerical approximations of these problems, such as the existence of a nonlinear potential term $f(u)$ and the presence of the  small interaction length $\epsilon$.
An appropriate numerical resolution of the Allen-Cahn equation requires a proper relation between physical and numerical scales, that is, the spatial size $h$ and the time size $k$ must be related to the perturbation parameter $\epsilon$.

In the past thirty years, there have been a large body of works on numerical simulations of  the Allen-Cahn equation (1.1)(cf.\cite{13,14,15,16,17,18} and the references therein). However, most of these works were conducted for a fixed parameter $\epsilon$. The error estimates,
which are deduced using the  Gronwall inequality \cite{Li, song}, depended on $\frac{1}{\epsilon}$ in exponential order. Indeed, such an estimate is obviously not useful for small parameter $\epsilon$, in particular, in discussing whether the flow of the computed numerical interfaces converge to the curvature-driven geometric flow. Less commonly investigated are error estimates which show an depend on $\frac{1}{\epsilon}$ in some (low) polynomial orders.
In general, the numerical analysis depending on $\frac{1}{\epsilon}$ in some (low) polynomial orders  can be significantly more difficult than that in exponential order. Nevertheless, such work has been reported in the following articles \cite{16,17,18,Feng_Li15}. One of the important ideas employed in the mentioned works is to use  a discrete spectrum estimate to derive  error order. In fact, the first such polynomial
order in $\frac{1}{\epsilon}$ a priori estimate was obtained by Feng and Prohl$^{\cite{16}}$ in 2003 for the finite element methods of the Allen-Cahn equation. And in 2015, Feng and Li$^{\cite{Feng_Li15}}$ developed fully discrete interior penalty discontinuous Galerkin methods for the Allen-Cahn equation, which is first-order-accurate-in-time numerical scheme and not unconditionally energy-stable scheme.
However, an essential feature of the Allen-Cahn equation is that it satisfies the energy laws (1.6), so it is important to design efficient and accurate numerical schemes that satisfy a corresponding discrete energy law, or in other words, energy stable.

In contrast to the papers referenced above, we propose a second-order-accurate-in-time,  unconditionally energy-stable with respect to the time and space step sizes, fully  discrete  interior penalty discontinuous Galerkin  finite element scheme for the Allen-Cahn problem (1.1)-(1.4). We  develop an interior penalty discontinuous Galerkin  finite element methods based on modified Crank-Nicolson scheme and a second-order-in-time approximation of the potential term $f(u)$, and establish polynomial order in $\frac{1}{\epsilon}$ a priori error estimates as well as to prove convergence and rates of convergence for the IPDGFE numerical interfaces. To the best of our knowledge,
no such numerical scheme and convergence  analysis for the Allen-Cahn problem (1.1)-(1.4) is available in the literature. The highlights of this paper include not only  presenting a second-order-accurate-in-time and unconditionally energy-stable scheme, but also
 using a discrete spectrum estimate to handle the midpoint of the discrete solutions $u^m$ and $u^{m+1}$ in the nonlinear term to achieve error bounds depending on  reciprocal of the perturbation parameter $\epsilon$  only in some lower polynomial order. Thus, it can be seen that the paper is not trivial extension of the article \cite{Feng_Li15} by Feng and Li.

The remainder of this paper is organized as follows. Section 2 includes a brief description of notions, and we recall a few facts and lemmas about the problem (1.1)-(1.4). In section 3, we present a fully discrete, nonlinear interior penalty discontinuous Galerkin method which is a second-order-in-time scheme based on a mid-point approximation of the potential term and proved to be unconditionally energy-stable and uniquely solvable, and provide a rigorous proof of convergence results for the proposed numerical method. In section 4, we prove the convergence and rates of convergence for the numerical interfaces of the numerical solutions to the sharp interface of the curvature-driven geometric flow. Finally, section 5 presents some of our numerical experiments to gauge the performance of the proposed interior penalty discontinuous Galerkin method.

\section{Preliminaries}\label{sec-2}
Let $\cT_h$ be a quasi-uniform ``triangulation" of $\Omega$ such that
$\overline{\Ome}=\bigcup_{K\in\cT_h} \overline{K}$. Let $h_K$ denote
the diameter of $K\in \cT_h$ and $h:=\mbox{max}\{h_K; K\in\cT_h\}$.
We recall that the standard broken Sobolev space $H^s(\cT_h)$ and DG finite
element space $V_h$ are defined as
\[
H^s(\mathcal{T}_h):=\prod_{K\in\cT_h} H^{s}(K), \qquad
V_h:=\prod_{K\in\cT_h} P_r(K),
\]
where $P_r(K)$ denotes the set of all polynomials whose degrees do not
exceed a given positive integer $r$.  Let $\mce^I$ denote the set of all
interior faces/edges of $\mct$, $\mce^B$ denote the set of all boundary
faces/edges of $\mct$, and $\mce:=\mce^I\cup \mce^B$. The $L^2$-inner product
for piecewise functions over the mesh $\cT_h$ is naturally defined by
\[
(u,v)_{\mathcal{T}_h}:= \sum_{K\in \cT_h} \int_{K} u v\, dx,
\]
and for any set $\mathcal{S}_h \subset \mce$, the $L^2$-inner product
over $\mathcal{S}_h$ is defined by
\begin{align*}
\bl u,v\br_{\mathcal{S}_h} :=\sum_{e\in \mathcal{S}_h} \int_e uv\, ds.
\end{align*}

Let $K, K'\in \cT_h$ and $e=\partial K\cap \partial K'$ and assume
global labeling number of $K$ is smaller than that of $K'$.
We choose $n_e:=n_K|_e=-n_{K'}|_e$ as the unit normal on $e$ and
define the following standard jump and average notations
across the face/edge $e$:
\begin{alignat*}{4}
[v] &:= v|_K-v|_{K'}
\quad &&\mbox{on } e\in \mce^I,\qquad
&&[v] :=v\quad
&&\mbox{on } e\in \mce^B,\\
\{v\} &:=\frac12\bigl( v|_K +v|_{K'} \bigr) \quad
&&\mbox{on } e\in \mce^I,\qquad
&&\{v\}:=v\quad
&&\mbox{on } e\in \mce^B
\end{alignat*}
for $v\in V_h$.
Let $M$ be a (large) positive integer. Define $\tau:=T/M$ and $t_m:=m\tau$
for $m=0,1,2,\cdots,M$ be a uniform partition of $[0,T]$. For a sequence
of functions $\{v^m\}_{m=0}^M$, we define the (backward) difference operator
\begin{equation}
d_t u^m:= \frac{u^m-u^{m-1}}{k}, \qquad m=1,2,\cdots,M. \nonumber
\end{equation}
First, we introduce the DG elliptic projection operator $P_r^h: H^s(\cT_h)\to V_h$ by
\begin{equation}\label{eq2.1}
a_h(v-P_r^h v, w_h) + \bigl( v-P_r^h v, w_h \bigr)_{\cT_h} =0
\quad\forall w_h\in V_h
\end{equation}
for any $v\in H^s(\cT_h)$.

We start with a well-known fact \cite{18} that the
Allen-Cahn equation \eqref{eq1.1} can be interpreted as the $L^2$-gradient
flow for the following Cahn-Hilliard energy functional
\begin{equation}\label{eq2.2}
J_\epsilon(v):= \int_\Omega \Bigl( \frac12 |\nabla v|^2
+ \frac{1}{\epsilon^2} F(v) \Bigr)\, dx.
\end{equation}

The following assumptions on the initial datum $u_0$ are made as in  \cite{ Feng_Li15, feng2014finite, feng2015analysis,feng2017finite,16,  li2015numerical,xu2016convex, 14} to derive a priori solution estimates.

\medskip
{\bf General Assumption} (GA)
\begin{itemize}
\item[(1)] There exists a nonnegative constant $\sigma_1$ such that
\begin{equation}\label{eq2.3}
J_{\epsilon}(u_0)\leq C\epsilon^{-2\sigma_1}.
\end{equation}
\item[(2)] There exists a nonnegative constant $\sigma_2$ such that
\begin{equation}\label{eq2.4}
\|\Delta u_0 -\epsilon^{-2} f(u_0)\|_{L^2(\Omega)} \leq C\epsilon^{-\sigma_2}.
\end{equation}

\item[(3)]
There exists nonnegative constant $\sigma_3$ such that
\begin{equation}\label{eq2.5}
\lim_{s\rightarrow0^{+}} \|\nabla u_t(s)\|_{L^2(\Omega)}\leq C\epsilon^{-\sigma_3}.
\end{equation}

\end{itemize}

The following solution estimates can be found in \cite{16, Feng_Li15}.

\begin{proposition}\label{prop2.1}
Suppose that \eqref{eq2.3} and \eqref{eq2.4} hold. Then the solution $u$ of
problem \eqref{eq1.1}--\eqref{eq1.4} satisfies the following estimates:
\begin{align} \label{eq2.6}
&\esssupI \|u(t)\|_{L^\infty(\Ome)} \leq 1,\\
&\esssupI\, J_{\epsilon}(u)
+ \int_{0}^{\infty} \|u_t(s)\|_{L^2(\Omega)}^2\, ds
\leq C \eps^{-2\sigma_1},\label{eq2.7} \\
&\int_{0}^{T} \|\Delta u(s)\|^2\, ds \leq C \eps^{-2(\sigma_1+1)}, \label{eq2.8}\\
&\esssupI \Bigl( \|u_t\|_{L^2(\Omega)}^2 +\|u\|_{H^2(\Omega)}^2 \Bigr)
+\int_{0}^{\infty} \|\nabla u_t(s)\|_{L^2(\Omega)}^2\, ds
\leq C \eps^{-2\max\{\sigma_1+1,\sigma_2\}}, \label{eq2.9} \\
&\int_{0}^{\infty} \Bigl(\|u_{tt}(s)\|_{H^{-1}(\Omega)}^2
+\|\Delta u_t(s)\|_{H^{-1}(\Omega)}^2\Bigr) \, ds
\leq C \eps^{-2\max\{\sigma_1+1,\sigma_2\}}. \label{eq2.10}
\end{align}
In addition to \eqref{eq2.3} and \eqref{eq2.4},
suppose that \eqref{eq2.5} holds, then $u$ also satisfies
\begin{align} \label{eq2.11}
&\esssupI \|\nabla u_t\|_{L^2(\Ome)}^2 +\int_0^{\infty} \|u_{tt}(s)\|_{L^2}^2 \,ds
\leq C\eps^{-2\max\{\sigma_1+2,\sigma_3\}},\\
&\int_{0}^{\infty} \|\Delta u_t(s)\|_{L^2(\Ome)}^2 \,ds
\leq C\eps^{-2\max\{\sigma_1+2, \sigma_3\}}. \label{eq2.12}
\end{align}

\end{proposition}

Next, we quote the following well known error estimate results from
$[21, 22]$.
\begin{lemma}\label{lem2.2}
Let $v\in W^{s,\infty}(\cT_h)$, then there hold
{\small
\begin{align}\label{eq2.13}
\|v-P_r^h v\|_{L^2(\cT_h)} +h\|\nab(v-P_r^h v)\|_{L^2(\cT_h)}
&\leq Ch^{\min\{r+1,s\}}\|u\|_{H^s(\cT_h)},\\
\frac{1}{|\ln h|^{\overline{r}}} \|v-P_r^h v\|_{L^\infty(\cT_h)}
+ h\|\nab(u-P_r^h u)\|_{L^\infty(\cT_h)}
&\leq Ch^{\min\{r+1,s\}}\|u\|_{W^{s,\infty}(\cT_h)}. \label{eq2.14}
\end{align}
}
where $\overline{r}:=\min\{1, r\}-\min\{1, r-1\}$.
\end{lemma}
\\
\\
Define $C_1$ as
\begin{equation}\label{eq2.15}
C_1=\underset{|\xi|\leq2}{\rm{max}}|f''(\xi)|.
\end{equation}
and $\hP_r^h$, corresponding to $P_r^h$, denote the elliptic projection
operator on the finite element space $S_h:=V_h\cap C^0(\overline{\Ome})$,
there holds the following estimate from $\cite{21}$:
\begin{equation}\label{eq2.16}
\|u-\hP_r^hu\|_{L^{\infty}}\leq Ch^{2-\frac{d}{2}}||u||_{H^2}.
\end{equation}

We now state our discrete spectrum estimate for the DG approximation.

\begin{proposition}\label{prop2.3}
Suppose there exists a positive number $\gamma>0$ such that the solution $u$ of
problem \eqref{eq1.1}--\eqref{eq1.4} satisfies
\begin{equation}\label{eq2.17}
\underset{t\in [0,T]}{\mbox{\rm ess sup}}\, \|u(t)\|_{W^{r+1,\infty}(\Ome)}
\leq C\eps^{-\gamma}.
\end{equation}
Then there exists an $\eps$-independent and $h$-independent constant
$c_0>0$ such that for $\eps\in(0,1)$ and a.e. $t\in [0,T]$
\begin{align}\label{eq2.18}
\lambda_h^{\mbox{\tiny DG}}(t):=\inf_{\psi_h\in V_h\atop\psi_h\not\equiv 0}
\frac{ a_h(\psi_h,\psi_h) + \frac{1}{\eps^2}\Bigl( f'\bigl(P_r^h u(t)\bigr)\psi_h,
\psi_h \Bigr)_{\cT_h}}{\|\psi_h\|_{L^2(\cT_h)}^2} \geq -c_0,
\end{align}
provided that $h$ satisfies the constraint
\begin{align}\label{eq2.19}
h^{2-\frac{d}{2}}
&\leq C_0 (C_1C_2)^{-1}\eps^{\max\{\sigma_1+3,\sigma_2+2\}},\\
&h^{\min\{r+1,s\}}|\ln h|^{\overline{r}} \leq C_0 (C_1C_2)^{-1}\eps^{\gamma+2},
\label{eq2.20}
\end{align}
where $C_2$ arises from the following inequality:
\begin{align}\label{eq2.21}
&\|u-P^h_r u\|_{L^{\infty}((0,T);L^{\infty}(\Ome)}
\leq C_2 h^{\min\{r+1,s\}}|\ln h|^{\overline{r}} \eps^{-\gamma},\\
&\|u-\hP^h_r u\|_{L^{\infty}((0,T);L^{\infty}(\Ome)}
\leq C_2 h^{2-\frac{d}{2}} \eps^{-\max\{\sigma_1+1,\sigma_2\}}. \label{eq2.22}
\end{align}
\end{proposition}
\begin{lemma}\label{lem2.4}
Let $\{S_{\ell} \}_{\ell\geq 1}$ be a positive nondecreasing sequence and
$\{b_{\ell}\}_{\ell\geq 1}$ and $\{k_{\ell}\}_{\ell\geq 1}$ be nonnegative sequences,
and $p>1$ be a constant. If
\begin{eqnarray}\label{eq2.23}
&S_{\ell+1}-S_{\ell}\leq b_{\ell}S_{\ell}+k_{\ell}S^p_{\ell} \qquad\mbox{for } \ell\geq 1,
\\ \label{eq2.24}
&S^{1-p}_{1}+(1-p)\mathop{\sum}\limits_{s=1}^{\ell-1}k_{s}a^{1-p}_{s+1}>0
\qquad\mbox{for } \ell\geq 2,
\end{eqnarray}
then
\begin{equation}\label{eq2.25}
S_{\ell}\leq \frac{1}{a_{\ell}} \Bigg\{S^{1-p}_{1}+(1-p)
\sum_{s=1}^{\ell-1}k_{s}a^{1-p}_{s+1}\Bigg\}^{\frac{1}{1-p}}\qquad\text{for}\ \ell\geq 2,
\end{equation}
where
\begin{equation}\label{eq2.26}
a_{\ell} := \prod_{s=1}^{\ell-1} \frac{1}{1+b_{s}} \qquad\mbox{for } \ell\geq 2.
\end{equation}
\end{lemma}


\section{Fully discrete IP-DG approximations}\label{sec-3}

\subsection{Discretized DG scheme} \label{sec-3.1}
We are now ready to introduce our fully discrete DG finite element methods
for problem \eqref{eq1.1}--\eqref{eq1.4}. They are defined by seeking
$u^m\in V_h$ for $m=0,1,2,\cdots, M$ such that
\begin{alignat}{2}\label{eq3.1}
\bigl( d_t u^{m+1},v_h\bigr)_{\mct}
+a_h(u^{m+\frac{1}{2}},v_h)+\frac{1}{\eps^2}\bigl( f^{m+1},v_h\bigr)_{\mct} &=0
&&\quad\forall v_h\in V_h,
\end{alignat}
where
\begin{align}\label{eq3.2}
a_h(u,v_h) &:=\bigl( \nabla u,\nabla v_h\bigr)_{\mct}
-\bigl\langle \{\p_n u\}, [v_h] \bigr\rangle_{\mce^I}  \\
&\hskip 1.1in
+\lambda \bigl\langle [u], \{\p_n v_h\} \bigr\rangle_{\mce^I} + j_h(u,v_h), \nonumber
\end{align}
\begin{align}\label{eq3.3}
j_h(u,v_h)&:=\sum_{e\in\mce^I}\frac{\sigma_e}{h_e}\bl [u],[v_h] \br_e,
\end{align}
\begin{align}\label{eq3.4}
f^{m+1}&:=\frac{1}{4} \bigl[(u^{m+1})^3+(u^{m+1})^2 u^m+u^{m+1}(u^{m})^2+(u^{m})^3\bigl]-u^{m+\frac{1}{2}}\\
&=\frac{F(u^{m+1})-F(u^m)}{u^{m+1}-u^m}.\nonumber
\end{align}
where  $u^{m+\frac{1}{2}}=\frac{u^{m+1}+u^{m}} {2}$, $\lambda=0,\pm1$ and $\sigma_e$ is a positive piecewise constant
function on $\mce^I$, which will be chosen later (see Lemma \ref{lem-3.1}).  In addition,
we need to supply $u_h^0$ to start the time-stepping,
whose choice will be clear (and will be specified) below.
\begin{lemma}\label{lem-3.1}
There exist constants $\sigma_0, \alpha>0$ such that for
$\sigma_e>\sigma_0$ for all $e\in \cE_h$ there holds
\begin{equation}\label{eq3.11a}
\Phi^h(v_h)\geq \alpha \|v_h\|_{1,\DG}^2 \qquad\forall v_h\in V_h, \nonumber
\end{equation}
where
\begin{equation}\label{eq3.11b}
\|v_h\|_{1,\DG}^2 := \|\nab v_h\|_{L^2(\cT_h)}^2 + j_h(v_h,v_h). \nonumber
\end{equation}
\end{lemma}

Now we introduce three mesh-dependent energy functionals
which can be regarded as DG counterparts of the continuous
Cahn-Hilliard energy $J_\eps$ defined in \eqref{eq2.2}.
\begin{align}\label{eq3.5}
\Phi^h(v) &:=\frac{1}{2} \|\nab v\|_{L^2(\cT_h)}^2
-\bigl\langle \{\p_n v\}, [v] \bigr\rangle_{\cE_h^I} + \frac12 j_h(v,v) \qquad
\forall v\in H^2(\cT_h), \\
J_\eps^h(v) &:= \Phi^h(v) +\frac{1}{\eps^2} \bigl( F(v), 1\bigr)_{\cT_h}
\qquad \forall v\in H^2(\cT_h), \label{eq3.6} \\
I_\eps^h(v) &:= \Phi^h(v) +\frac{1}{\eps^2} \bigl( F_c^+(v), 1\bigr)_{\cT_h}
\qquad \forall v\in H^2(\cT_h), \label{eq3.7}
\end{align}
It is easy to check that $\Phi^h$ and $I_\eps^h$ are convex functionals
but $J_\eps^h$ is not because $F$ is not convex. Moreover, we have:
\begin{lemma}\label{lem3.2}
Let $\lambda =-1$ in \eqref{eq3.2}, then there holds for all $v_h,w_h\in V_h$
\begin{align}\label{eq3.8}
\Bigl(\frac{\delta \Phi^h(v_h)}{ \delta v_h}, w_h \Bigr)_{\cT_h}
&:=\lim_{s\to 0} \frac{\Phi^h(v_h+ s w_h)-\Phi^h(v_h) }{s}
=a_h(v_h, w_h), \\
\Bigl(\frac{\delta J_\eps^h(v_h)}{ \delta v_h}, w_h \Bigr)_{\cT_h}
:&=\lim_{s\to 0} \frac{J_\eps^h(v_h+ s w_h)-J_\eps^h(v_h) }{s} \label{eq3.9} \\
&=a_h(v_h, w_h) +\frac{1}{\eps^2} \bigl(F^{\prime}(v_h), w_h \bigr)_{\cT_h},
\nonumber\\
\Bigl(\frac{\delta I_\eps^h(v_h)}{ \delta v_h}, w_h \Bigr)_{\cT_h}
:&=\lim_{s\to 0} \frac{I_\eps^h(v_h+ s w_h)-I_\eps^h(v_h) }{s} \\
&=a_h(v_h, w_h) +\frac{1}{\eps^2} \bigl((F_c^+)^{\prime}(v_h), w_h \bigr)_{\cT_h}.
\nonumber
\end{align}
\end{lemma}
\subsection{Stability of the DG scheme} \label{sec-3.2}
\begin{theorem}
The scheme \eqref{eq3.1}--\eqref{eq3.4} is unconditionally stable for all $h,k>0$ .
\end{theorem}
Proof: We have the DG scheme as below:
\begin{align}\label{eq3.11}
\big(d_tu^{m+1},v_h\bigr)+a_h\big(\frac{u^{m+1}+u^{m}}{2},v_h\big)+\frac{1}{\epsilon^2}\big(f^{m+1},v_h\big)=0.
\end{align}
Let $v_h=d_tu^{m+1}$, and we will get:
\begin{align}\label{eq3.12}
\big(d_tu^{m+1},d_tu^{n+1}\bigr)+a_h\big(\frac{u^{m+1}+u^{m}}{2},d_tu^{m+1}\big)&\\
+\frac{1}{\epsilon^2}\big(\frac{F(u^{m+1})-F(u^m)}{u^{m+1}-u^m},d_tu^{m+1}\big)&=0.\notag
\end{align}
Rearrange it to get:
\begin{align}\label{eq3.13}
\|d_tu^{m+1}\|^2_{L^2}+\frac 1 2 d_t[a_h(u^{m+1},u^{m+1})]+\frac{1}{\epsilon^2}d_t F(u^{m+1})=0,
\end{align}
\begin{align}\label{eq3.14}
d_t[ \frac1 2a_h(u^{m+1},u^{m+1})+\frac{1}{\epsilon^2}(F(u^{m+1}),1)]\leq0.
\end{align}
And the proof is complete.
\subsection{Well-posedness of the DG scheme}\label{sec-3.3}
We want to get a second order approximation of $f(u^{m+1}, u^m)$, which leads to unconditionally energy stable schemes. We split the function $F(v)=\frac14 (v^2-1)^2$ into the difference of two convex parts and get the convex decomposition $F(v)=F_c^+(v)-F_c^-(v)$,where $F_c^+(v):= \frac14 (v^4+1)$and $F_c^-(v):=\frac12 v^2$.\\
\\
Now we want to construct a second-order energy-stable scheme to approximate the two convex functions $F_c^+(u)$ and $F_c^-(u)$.
\begin{align}
f^+(u^{m+1},u^m)=\frac{F_c^+(u^{m+1})-F_c^+(u^m)}{u^{m+1}-u^m},\nonumber
\end{align}
\begin{align}
f^-(u^{m+1},u^m)=\frac{F_c^-(u^{m+1})-F_c^-(u^m)}{u^{m+1}-u^m}.\nonumber
\end{align}
\begin{theorem}
Under the constraint $k<2\epsilon^2$, there exists a unique solution of the scheme \eqref{eq3.1}- \eqref{eq3.4}.
\end{theorem}

Proof: Define the following functional:
\begin{align}\label{eq3.15}
&J(u^{m+1})=\frac1 4 a_h(u^{m+1},u^{m+1})+ \frac{1}{\epsilon^2}\int_{\cT_h}F_+(u^{m+1},u^m)\\
&+(\frac{1}{2k}-\frac{1}{4\epsilon^2})\|u^{m+1}\|^2_{L^2(\cT_h)}+\frac1 2 a_h(u^m,u^{m+1})+\int_{\cT_h}(-\frac{1}{2\epsilon^2}-\frac{1}{k})u^mu^{m+1}.\nonumber
\end{align}
Take the  derivative of the functional $J(u^{m+1})$, and will get:
\begin{align}\label{eq3.16}
&\Bigl(\frac{\delta J(u^{m+1})}{ \delta u^{m+1}}, v_h \Bigr)_{\cT_h}=\frac1 2 a_h(u^{m+1},v_h)+\frac{1}{\epsilon^2}\int_{\cT_h}f^+(u^{m+1},u^m)\\
&+(\frac{1}{2k}-\frac{1}{4\epsilon^2})2(u^{m+1},v_h)_{\cT_h}+\frac 1 2 a_h(u^m,v_h)+(-\frac{1}{2\epsilon^2}-\frac{1}{k})(u^{m},v_h)_{\cT_h}.\nonumber
\end{align}
Rearrange it, and we will get:
\begin{align}\label{eq3.17}
\Bigl(\frac{\delta J(u^{m+1})}{ \delta u^{m+1}}, v_h \Bigr)_{\cT_h}=\bigl( d_t u^{m+1},v_h\bigr)_{\mct}
+a_h(u^{m+\frac{1}{2}},v_h)+\frac{1}{\eps^2}\bigl( f^{m+1},v_h\bigr)_{\mct} &=0.
\end{align}
Also we can see from $(3.10)$ the first two terms of $J(u^{m+1})$ are convex, also since the last two terms are linear with respect to $u^{m+1}$, so they are also convex, so if we restrict the coefficient of third term to be positive, that is, if we restrict $k<2\epsilon^2$, then the $J(u^{m+1})$ will be a convex functional, and the uniqueness of the solution to this scheme is approved.

\subsection{Error estimates analysis}\label{sec-3.4}
The main result of this subsection is the following error
estimate theorem.

\begin{theorem}\label{thm3.5}
suppose $\sigma_e>\max\{\sigma_0,\sigma_0^{\prime}\}$.
Let $u$ and $\{u_h^m\}_{m=1}^M$ denote respectively the solutions of problems
\eqref{eq1.1}--\eqref{eq1.4} and \eqref{eq3.1}--\eqref{eq3.5}.
Assume $u\in H^2((0,T);$ $L^2(\Ome))\cap L^2((0,T); W^{s,\infty}(\Ome))$
and suppose (GA) and \eqref{eq2.17} hold. Then, under the following
mesh and initial value constraints:

\begin{align*}
h^{2-\frac{d}{2}} &\leq C_0 (C_1C_2)^{-1}\eps^{\max\{\sigma_1+3,\sigma_2+2\}},
\end{align*}
\begin{align*}
h^{\min\{r+1,s\}}|\ln h|^{\overline{r}} &\leq C_0 (C_1C_2)^{-1}\eps^{\gamma+2},
\end{align*}
\begin{align*}
k<A(\epsilon),
\end{align*}
\begin{align*}
u_h^0\in S_h\mbox{ such that }\quad
\|u_0 -u_h^0\|_{L^2(\cT_h)} &\leq C h^{\min\{r+1,s\}},
\end{align*}
there hold
\begin{align}\label{eq3.18}
\max_{0\leq m\leq M} \|u(t_m)-u_h^m\|_{L^2(\cT_h)}
&\leq C(k^2+h^{\min\{r+1,s\}})\eps^{-(\sigma_1+2)}.\\
\Bigl( k\sum_{m=1}^M \| u(t_m)-u_h^m\|_{H^1(\cT_h)}^2 \Bigr)^{\frac{1}{2}}
&\leq C(k^2+h^{\min\{r+1,s\}-1})\eps^{-(\sigma_1+3)}, \label{eq3.19} \\
\max_{0\leq m\leq M} \|u(t_m)-u_h^m\|_{L^\infty(\cT_h)}
&\leq C h^{\min\{r+1,s\}} |\ln h|^{\overline{r}} \eps^{-\gamma} \label{eq3.20} \\
&\qquad
+Ch^{-\frac{d}{2}}(k^2+h^{\min\{r+1,s\}})\eps^{-(\sigma_1+2)}.  \nonumber
\end{align}
\end{theorem}


Proof: Since the proof is long, we split the proof into four steps:\\
\medskip
{\em Step 1}: \\
We write:
\[
u(t_m)-u^m=\eta^m + \xi^m,\quad \eta^m:=u(t_m)-P_r^h u(t_m),\quad
\xi^m:=P_r^h u(t_m)- u^m.
\]
Multiply  $v_h$ on both sides of the Allen-Cahn equation in $(1.1)$ at the point $u(t_{m+\frac{1}{2}})$
\begin{align}\label{eq3.21}
\bigl(u_t(t_{m+\frac{1}{2}}),v_h\bigr)_{\cT_h} +a_h(u(t_{m+\frac{1}{2}}),v_h)
+\frac{1}{\eps^2}\bigl( f(u(t_{m+\frac{1}{2}})), v_h\bigr)_{\cT_h}
= 0,
\end{align}
for all $ v_h\in V_h$, where $t_{m+\frac{1}{2}}=\frac{t_{m+1}+t_{m}} {2}$.\\
\\
Subtract $(3.1)$ from $(3.21)$, we get the following equation:
\begin{align}\label{eq3.22}
&\bigl(u_t(t_{m+\frac{1}{2}})- \frac{u^{m+1}-u^{m}}{k},v_h\bigr)_{\cT_h}
+a_h\big(u(t_{m+\frac{1}{2}})-\frac{u^{m+1}+u^{m}} {2},v_h\big)\\
\hskip 0.0in
&+\frac{1}{\eps^2}\bigl( f(u(t_{m+\frac{1}{2}}))-f^{m+1}, v_h \bigr)_{\cT_h} = 0.\nonumber
\end{align}
From Taylor expansion:
\[
u(t_{m+1})=u(\frac{t_{m+1}+t_{m}}{2})+u_t(\frac{t_{m+1}+t_{m}}{2})(\frac{t_{m+1}-t_{m}}{2})+R_1^{m},
\]
where $R_1^{m}=u_{tt}(\xi_1)(\frac{t_{m+1}-t_{m}}{2})^2$.
\[
u(t_{m})=u(\frac{t_{m+1}+t_{m}}{2})-u_t(\frac{t_{m+1}+t_{m}}{2})(\frac{t_{m+1}-t_{m}}{2})+{R_2^{m}},
\]
where $R_2^{m}=u_{tt}(\xi_2)(\frac{t_{m+1}-t_{m}}{2})^2$.
And we will get:\\
\begin{align}\label{eq3.23}
u(t_{m+\frac{1}{2}})=\frac{u(t_{m+1})+u(t_{m})}{2}-\frac{\big(R_1^{m}+R_2^{m}\big)}{2},
\end{align}
\begin{align}\label{eq3.24}
u_t(t_{m+\frac{1}{2}})=\frac{u(t_{m+1})-u(t_{m})}{k}-\frac{\big(R_1^{m}-R_2^{m}\big)}{k}.
\end{align}
Use $(3.23)$ and $(3.24)$  into $(3.22)$ , we will get:
\begin{align}\label{eq3.25}
&\bigl(\frac{\xi^{m+1}-\xi^{m}}{k}+\frac{\eta^{m+1}-\eta^{m}}{k}-\frac{\big(R_1^{m}-R_2^{m}\big)}{k},v_h\bigr)_{\cT_h} \\
\hskip 0.0in
&+a_h\bigl(\frac{\xi^{m+1}+\xi^{m}}{2}+\frac{\eta^{m+1}+\eta^{m}}{2}-\frac{\big(R_1^{m}+R_2^{m}\big)}{2},v_h\bigr)\nonumber \\
\hskip 0.0in
&+\frac{1}{\eps^2}\bigl( f(u(t_{m+\frac{1}{2}}))-f^{m+1}, v_h \bigr)_{\cT_h} = 0.\nonumber
\end{align}
\begin{align}\label{eq3.26}
&\bigl(d_t \xi^{m+1},v_h\bigr)_{\cT_h} +a_h(\frac{\xi^{m+1}+\xi^{m}}{2},v_h)\\
&+\frac{1}{\eps^2}\bigl( f(u(t_{m+\frac1 2}))-f^{m+1}, v_h\bigr)_{\cT_h} \nonumber\\
\hskip 0.5in
&= \bigl(\frac{\big(R_1^{m}-R_2^{m}\big)}{k},v_h\bigr)_{\cT_h}
-\bigl(d_t \eta^{m+1},v_h\bigr)_{\cT_h}\nonumber\\
&-a_h(\frac{\eta^{m+1}+\eta^{m}}{2},v_h)+a_h(\frac{\big(R_1^{m}+R_2^{m}\big)}{2},v_h) \nonumber\\
&= \bigl(\frac{\big(R_1^{m}-R_2^{m}\big)}{k},v_h\bigr)_{\cT_h}
-\bigl(d_t \eta^{m+1},v_h\bigr)_{\cT_h} \nonumber\\
&+(\frac{\eta^{m+1}+\eta^{m}}{2},v_h)_{\cT_h}+a_h(\frac{\big(R_1^{m}+R_2^{m}\big)}{2},v_h). \nonumber
\end{align}
Let $v_h=\frac{\xi^{m+1}+\xi^{m}}{2}$, for the first term on the left hand side:
\begin{align}\label{eq3.27}
&\bigl(d_t \xi^{m+1},\frac{\xi^{m+1}+\xi^{m}}{2}\bigr)_{\cT_h}= \frac12 d_t \|\xi^{m+1}\|_{L^2(\cT_h)}^2.
\end{align}
We split  the third term on the left hand side in $(3.26)$ into two parts and deal with them separately:
\begin{align}\label{eq3.28}
&\frac{1}{\eps^2}\bigl( f(u(t_{m+\frac1 2}))-f^{m+1}, \frac{\xi^{m+1}+\xi^{m}}{2}\bigr)_{\cT_h}\\
&=\frac{1}{\eps^2}\bigl( f(u(t_{m+\frac1 2}))-f(\frac{u(t_{m+1})+u(t_{m})}{2}), \frac{\xi^{m+1}+\xi^{m}}{2}\bigr)_{\cT_h} \nonumber\\
&+\frac{1}{\eps^2}\bigl( f(\frac{u(t_{m+1})+u(t_{m})}{2})-f^{m+1}, \frac{\xi^{m+1}+\xi^{m}}{2}\bigr)_{\cT_h}.\nonumber
\end{align}
Let $\hat{u}(t_{m+\frac{1}{2}})=\frac{u(t_{m+1})+u(t_{m})}{2}$, the we have the following:\\
\begin{align}\label{eq3.29}
 &f(u(t_{m+\frac1 2}))-f(\hat{u}(t_{m+\frac{1}{2}}))\\
 &=f\big(\hat{u}(t_{m+\frac{1}{2}})-\frac1 8k^2(u''(\xi_1)+u''(\xi_2))\big)-f(\hat{u}(t_{m+\frac{1}{2}}))\nonumber\\
 &=f'(\xi_{12})(-\frac1 8)k^2((u''(\xi_1)+u''(\xi_2))\geq-Ck^2.\nonumber
\end{align}
Since $f'$ and $u''$ both are bounded, we will get the following inequality by Cauchy-Schwarz inequality:
\begin{align}\label{eq3.30}
&\frac{1}{\epsilon^2}\bigl(f(u(t_{m+\frac1 2}))-f(\hat{u}(t_{m+\frac{1}{2}})),\frac{\xi^{m+1}+\xi^{m}}{2}\bigr)_{\cT_h}\\
&\geq -\frac{1}{\epsilon^2}\bigl(Ck^2,\frac{\xi^{m+1}+\xi^{m}}{2}\bigr)_{\cT_h}\nonumber\\
&\geq -\frac{1}{\epsilon^4}Ck^4-\|\xi^{m+\frac 1 2}\|^2_{L^2(\cT_h)}.\nonumber
\end{align}
For the last term of the right hand side in \eqref{eq3.26}:
\begin{align}\label{eq3.31}
&a_h\big(\frac{\big(R_1^{m}+R_2^{m}\big)}{2},\frac{\xi^{m+1}+\xi^{m}}{2}\big)=a_h\big(\frac{\big(R_1^{m}+R_2^{m}\big)}{2\epsilon},\frac{\epsilon(\xi^{m+1}+\xi^{m})}{2}\big)\\
&\leq a_h\big(\frac{(R_1^{m}+R_2^{m})}{2\epsilon},\frac{(R_1^{m}+R_2^{m})}{2\epsilon}\big)+a_h\big(\frac{\epsilon(\xi^{m+1}+\xi^{m})}{2},\frac{\epsilon(\xi^{m+1}+\xi^{m})}{2}\big)\nonumber\\
&\leq Ck^4\epsilon^{-2}+\epsilon^2a_h\big(\xi^{m+\frac1 2},\xi^{m+\frac1 2}\big).\nonumber
\end{align}
Substitute \eqref{eq3.27},\eqref{eq3.30} and \eqref{eq3.31} into \eqref{eq3.26}, and we will get:
\begin{align}\label{eq3.32}
&\frac12d_t \|\xi^{m+1}\|_{L^2(\cT_h)}^2 +a_h(\frac{\xi^{m+1}+\xi^{m}}{2},\frac{\xi^{m+1}+\xi^{m}}{2})\\
&+\frac{1}{\eps^2}\bigl( f(\hat{u}(t_{m+\frac{1}{2}}))-f^{m+1}, \frac{\xi^{m+1}+\xi^{m}}{2}\bigr)_{\cT_h}\nonumber\\
&= \bigl(\frac{\big(R_1^{m}-R_2^{m}\big)}{k},v_h\bigr)_{\cT_h} -\bigl(d_t \eta^{m+1},v_h\bigr)_{\cT_h} \nonumber\\
&+(\frac{\eta^{m+1}+\eta^{m}}{2},v_h)_{\cT_h}+a_h(\frac{R_1^{m}+R_2^{m}}{2},v_h) \nonumber\\
&\leq\bigl( \|\bigl(\frac{\big(R_1^{m}-R_2^{m}\big)}{k})\|_{L^2(\cT_h)}^2+\|d_t\eta^{m+1}\|_{L^2(\cT_h)}^2\nonumber\\
&+\|(\frac{\eta^{m+1}+\eta^{m}}{2})\|_{L^2(\cT_h)}^2\bigr)\bigl(\frac{\xi^{m+1}+\xi^{m}}{2}\big)\|_{L^2(\cT_h)}^2\nonumber\\
&+Ck^4[\epsilon^{-4}+\epsilon^{-2}]+\epsilon^2a_h\big(\xi^{m+\frac1 2},\xi^{m+\frac1 2}\big) +\|\xi^{m+\frac 1 2}\|^2_{L^2(\cT_h)}.  \nonumber
\end{align}
Using the integral form of Taylor formula, we can get:
\[
|\frac{R_1^{m}-R_2^{m}}{k}|=|\frac{k(u_{tt}(\xi_1)-u_{tt}(\xi_2))}{4}=|\frac{ku_{ttt}(\xi_{11})(\xi_1-\xi_2)}{4}|\leq Ck^2.
\]
Hence
\begin{align}\label{eq3.33}
\|\frac{R_1^{m}-R_2^{m}}{k}\|^2_{L^2(\cT_h)}\leq Ck^4.
\end{align}
Summing in m from 1 to $\ell$,  using \eqref{eq3.13},\eqref{eq3.32} and \eqref{eq3.33}, and we will get the following inequality:
\begin{align}\label{eq3.34}
&\|\xi^{\ell}\|_{L^2(\cT_h)}^2 + 2k\sum_{m=1}^\ell a_h(\frac{\xi^{m}+\xi^{m-1}}{2},\frac{\xi^{m}+\xi^{m-1}}{2})\\
&+2k\sum_{m=1}^\ell \frac{1}{\eps^2}\bigl( f(\hat{u}(t_{m-\frac{1}{2}}))-f^{m}, \frac{\xi^{m}+\xi^{m-1}}{2}\bigr)_{\cT_h}\nonumber\\
&\leq \|\xi^{0}\|_{L^2(\cT_h)}^2 +Ch^{2\min\{r+1,s\}}\, \|u\|_{H^1((0,T);H^s(\Ome))}^2  \nonumber\\
&+2Ck^4[\epsilon^{-4}+\epsilon^{-2}+1]+2k\sum_{m=1}^\ell \epsilon^2a_h\big(\xi^{m-\frac1 2},\xi^{m-\frac1 2}\big) +4k\sum_{m=1}^\ell \|\xi^{m-\frac1 2}\|^2_{L^2(\cT_h)}. \nonumber
\end{align}
\medskip
{\em Step 2}:
We want to bound the term$\bigl( f(\hat{u}(t_{m+\frac{1}{2}}))-f^{m+1}, \frac{\xi^{m+1}+\xi^{m}}{2}\bigr)_{\cT_h}$ on the left hand side of \eqref{eq3.34}:
\begin{align}\label{eq3.35}
f(\hat{u}(t_{m+\frac{1}{2}}))-f^{m+1}=&[f(\hat{u}(t_{m+\frac{1}{2}}))-f(P_r^h \hat{u}(t_{m+\frac{1}{2}}))]\\
 &+ [f(P_r^h \hat{u}(t_{m+\frac{1}{2}}))-f^{m+1}].\notag
\end{align}
For the first part on the right hand side of $(3.35)$,we get:
\begin{align}\label{eq3.36}
|f(u^{m+\frac 1 2})-f(P_r^h \hat{u}(t_{m+\frac{1}{2}}))| =&|f'(\xi)|\bigl|\hat{u}(t_{m+\frac{1}{2}})-P_r^h \hat{u}(t_{m+\frac{1}{2}})\bigr|\notag\\
&\geq -C|\frac{\eta^{m+1}+\eta^{m}}{2}|.
 \end{align}
 For the second part on the right hand side of \eqref{eq3.35},we get:
 \begin{align}\label{eq3.37}
&f(P_r^h \hat{u}(t_{m+\frac1 2}))-f^{m+1}\\
&=\bigl(\frac{P_r^h u(t_{m+1})+P_r^h u(t_{m})}{2}\bigr)^3-\bigl(\frac{P_r^h u(t_{m+1})+P_r^h u(t_{m})}{2}\bigr)\notag\\
 &-\bigl[\frac1 4[(u^{m+1})^3+(u^{m+1})^2u^m+u^{m+1}(u^m)^2+(u^m)^3]-\frac{u^{m+1}+u^m}{2}\bigr] \nonumber\\
&=\frac{\bigl(P_r^h u(t_{m+1})+P_r^h u(t_{m})\bigr)^3}{8}\notag\\
&-\frac2 8[(u^{m+1})^3+(u^{m+1})^2u^m+u^{m+1}(u^m)^2+(u^m)^3]\nonumber\\
&-\bigl[\bigl(\frac{P_r^h u(t_{m+1})+P_r^h u(t_{m})}{2}\bigr)-\frac{u^{m+1}+u^m}{2}\bigr] \nonumber\\
&=\frac{\bigl(P_r^h u(t_{m+1})+P_r^h u(t_{m})\bigr)^3}{8}-\frac2 8[(P_r^h u(t_{m+1})-\xi^{m+1})^3\notag\\
&+(P_r^h u(t_{m+1})-\xi^{m+1})^2(P_r^h u(t_{m})-\xi^{m})-\frac{(\xi^{(m+1)}+\xi^m)}{2}\nonumber\\
&+(P_r^h u(t_{m+1})-\xi^{m+1})(P_r^h u(t_{m})-\xi^{m})^2+(P_r^h u(t_{m})-\xi^{m})^3].\nonumber
 \end{align}
 We split the above into four terms: constant term with resect to $\xi^{m+1}$ and\ \ $\xi^{m}$, linear, quadratic and cubic in terms of $\xi^{m+1}$ and\  \    $\xi^{m}$.\\
 For constant term, we have
 \begin{align}\label{eq3.38}
&\bigl(\frac1 8 (P_r^h u(t_{m+1})-P_r^h u(t_{m}))^2(P_r^h u(t_{m+1})+P_r^h u(t_{m})), \frac{\xi^{m+1}+\xi^{m}}{2}\bigr)_{\cT_h}\\
&\geq -C(h^4+k^2)\bigl(1, \frac{\xi^{m+1}+\xi^{m}}{2}\bigr)_{\cT_h}\nonumber\\
&\geq -C(h^8+k^4)-C\|\xi^{m+\frac1 2}\|^2_{L^2(\cT_h)}.\nonumber
 \end{align}
 By the boundness of $P_r^h u^{m}$ and $|P_r^h u(t_{m+1})-P_r^h u(t_{m})|\leq h^2+k$.\\
 \\
 For the linear term, we have the following:\\
 \begin{align}\label{eq3.39}
& l= \frac 1 4\bigl\{ \xi^{m+1}[3(P_r^h u(t_{m+1}))^2+P_r^h u(t_{m+1})P_r^h u(t_{m})+(P_r^h u(t_{m}))^2]\\
&+\xi^{m}[3(P_r^h u(t_{m}))^2+P_r^h u(t_{m+1})P_r^h u(t_{m})+(P_r^h u(t_{m+1}))^2]\bigr\}-\frac{(\xi^{m+1}+\xi^m)}{2}\nonumber\\
& =\frac1 4(\xi^{m+1}+\xi^{m})(P_r^h u(t_{m+1})+P_r^h u(t_{m}))^2\nonumber\\
&+\frac1 2[\xi^{m+1}(P_r^h u(t_{m+1}))^2+\xi^{m}(P_r^h u(t_{m}))^2]-\frac{(\xi^{m+1}+\xi^m)}{2}.\nonumber
 \end{align}
  And we have:\\
 \begin{align}\label{eq3.40}
 &\bigl( \frac1 4(\xi^{m+1}+\xi^{m})(P_r^h u(t_{m+1})+P_r^h u(t_{m}))^2\\
 &+\frac 12[\xi^{m+1}(P_r^h u(t_{m+1}))^2+\xi^{m}(P_r^h u(t_{m}))^2], \frac{\xi^{m+1}+\xi^{m}}{2}\bigr)_{\cT_h}\nonumber\\
 &=\bigl( \frac 1 2( P_r^h u(t_{m+1})+P_r^h u(t_{m}))^2, (\frac{\xi^{m+1}+\xi^{m}}{2})^2\bigr)_{\cT_h}\nonumber\\
 &+(\frac 12[\xi^{m+1}(P_r^h u(t_{m+1}))^2+\xi^{m}(P_r^h u(t_{m}))^2], \frac{\xi^{m+1}+\xi^{m}}{2}\bigr)_{\cT_h}.\nonumber
 \end{align}
 By using the Schwarz Inequality and $|P_r^h u(t_{m+1})-P_r^h u(t_{m})|\leq C(h^2+k)$, we get the following inequalities for the first and second terms of the right hand side of \eqref{eq3.40}:
 \begin{align}\label{eq3.41}
 &\bigl( \frac 1 2( P_r^h u(t_{m+1})+P_r^h u(t_{m}))^2, (\frac{\xi^{m+1}+\xi^{m}}{2})^2\bigr)_{\cT_h}\\
&\geq\bigl(2(P_r^h u(t_{m}))^2, (\frac{\xi^{m+1}+\xi^{m}}{2})^2\bigr)_{\cT_h}-C(h^2+k)\|\xi^{m+\frac1 2}\|_{L^2(\cT_h)}^2.\nonumber
 \end{align}
  \begin{align}\label{eq3.42}
 &(\frac 12[\xi^{m+1}(P_r^h u(t_{m+1}))^2+\xi^{m}(P_r^h u(t_{m}))^2], \frac{\xi^{m+1}+\xi^{m}}{2}\bigr)_{\cT_h}\\
&\geq\bigl((P_r^h u(t_{m}))^2, (\frac{\xi^{m+1}+\xi^{m}}{2})^2\bigr)_{\cT_h}-C(h^2+k)\|\xi^{m+\frac 1 2}\|_{L^2(\cT_h)}^2.\nonumber
 \end{align}
   \begin{align}\label{eq3.43}
&\bigl(l, \frac{\xi^{m+1}+\xi^{m}}{2}\bigr)_{\cT_h}\\
&\geq\bigl(3(P_r^h u(t_{m}))^2-1, (\frac{\xi^{m+1}+\xi^{m}}{2})^2\bigr)_{\cT_h}-C(h^2+k)\|\xi^{m+\frac1 2}\|_{L^2(\cT_h)}^2\nonumber\\
&=\bigl((f'(P_r^h u(t_{m})), (\frac{\xi^{m+1}+\xi^{m}}{2})^2\bigr)_{\cT_h}-C(h^2+k)\|\xi^{m+\frac1 2}\|_{L^2(\cT_h)}^2.\nonumber
 \end{align}
 For the quadratic term,  we get the inequality below;\\
 \begin{align}\label{eq3.44}
 &q=3(\xi^{m+1})^2P_r^h u(t_{m+1})+(\xi^{m+1})^2P_r^h u(t_{m})+2\xi^{m+1}\xi^{m}P_r^h u(t_{m+1})\\
 &+(\xi^{m})^2P_r^h u(t_{m+1})+2\xi^{m+1}\xi^{m}P_r^h u(t_{m})+3(\xi^{m})^2P_r^h u(t_{m})\nonumber\\ &\geq -C_1[(\xi^{m+1})^2+(\xi^{m})^2].\nonumber
 \end{align}
 So we get:\\
\begin{align}\label{eq3.45}
&\bigl(q, \frac{\xi^{m+1}+\xi^{m}}{2}\bigr)_{\cT_h}\\
&\geq -C_1\bigl((\xi^{m+1})^2+(\xi^{m})^2, \frac{\xi^{m+1}+\xi^{m}}{2}\bigr)_{\cT_h}\nonumber\\
&\geq-C \|\xi^{m+\frac{1}{2}}\|_{L^3(\cT_h)}^3.\nonumber
 \end{align}
For cubic term, we have:\\
 \begin{align}\label{eq3.46}
&c= \frac1 4\bigl[(\xi^{m+1})^3+(\xi^{m+1})^2\xi^{m}+\xi^{m+1}(\xi^{m})^2+(\xi^{m})^3\bigr]\nonumber\\
&\    =\frac1 4\bigl[(\xi^{m+1})^2+(\xi^{m})^2\bigr](\xi^{m+1}+\xi^{m}),
 \end{align}
 Then we have
\begin{align}\label{eq3.47}
\bigl(c, \frac{\xi^{m+1}+\xi^{m}}{2}\bigr)_{\cT_h}=\bigl((\xi^{m+1})^2+(\xi^{m})^2, (\xi^{m+1}+\xi^{m})^2\bigr)_{\cT_h}\geq 0.
\end{align}
Combine all above together, we will get:
\begin{align}\label{eq3.48}
&\bigl(f(P_r^h \hat{u}(t_{m+\frac1 2}))-f^{m+1}, \frac{\xi^{m+1}+\xi^{m}}{2}\bigr)_{\cT_h}\\
&\geq -C|(\eta^{m+\frac1 2}, \xi^{m+\frac1 2})|_{\cT_h}-C(h^8+k^4)-C\|\xi^{m+\frac1 2}\|_{L^2(\cT_h)}\nonumber\\
&\bigl((f'(P_r^h u(t_{m}))), (\frac{\xi^{m+1}+\xi^{m}}{2})^2\bigr)_{\cT_h}-C(h^2+k)\|\xi^{m+\frac1 2}\|_{L^2(\cT_h)}^2\nonumber\\
&-C \|\xi^{m+\frac{1}{2}}\|_{L^3(\cT_h)}^3+\frac{4k}{\epsilon^2}\bigl((\xi^{m+1})^2+(\xi^{m})^2, (\xi^{m+1}+\xi^{m})^2\bigr)_{\cT_h}.\nonumber
\end{align}
Summing in m a from 1 to $\ell$ and we will get the following:
\begin{align}\label{eq3.49}
&\frac{2k}{\epsilon^2}\sum_{m=1}^\ell \bigl( f(P_r^h \hat{u}(t_{m+\frac1 2}))-f^{m+1}, \frac{\xi^{m+1}+\xi^{m}}{2}\bigr)_{\cT_h}\\
&\geq -\frac{Ck}{\epsilon^2}\sum_{m=1}^\ell \|\eta^{m+\frac1 2}\|_{\cT_h}\|\xi^{m+\frac1 2}\|_{\cT_h}-C\frac{1}{\epsilon^2}(h^8+k^4)-\frac{Ck}{\epsilon^2}\sum_{m=1}^\ell \|\xi^{m}\|^2_{L^2(\cT_h)}\nonumber\\
&+\frac{2k}{\epsilon^2}\sum_{m=1}^\ell \bigl(f'(P_r^h u(t_{m})), (\frac{\xi^{m+1}+\xi^{m}}{2})^2\bigr)_{\cT_h}-C\frac{k}{\epsilon^2}  (h^2+k)\sum_{m=1}^\ell\|\xi^{m}\|_{L^2(\cT_h)}^2\nonumber\\
&-C\frac{k}{\epsilon^2} \sum_{m=1}^\ell \|\xi^{m+\frac{1}{2}}\|_{L^3(\cT_h)}^3
+\sum_{m=1}^\ell \bigl((\xi^{m})^2+(\xi^{m+1})^2, (\xi^{m}+\xi^{m+1})^2\bigr)_{\cT_h},\nonumber\\
&\geq-C h^{2\min\{r+1,s\}} \eps^{-4} \|u\|_{L^2((0,T);H^s(\Ome)}^2-C\frac{1}{\epsilon^2}(h^8+k^4)\nonumber\\
& +\frac{2k}{\epsilon^2}\sum_{m=1}^\ell \bigl(f'(P_r^h u(t_{m})), (\frac{\xi^{m+1}+\xi^{m}}{2})^2\bigr)_{\cT_h}\nonumber\\
&+\frac{2k}{\epsilon^2}\sum_{m=1}^\ell \bigl((\xi^{m})^2+(\xi^{m+1})^2, (\xi^{m}+\xi^{m+1})^2\bigr)_{\cT_h}-C\frac{k}{\epsilon^2} \sum_{m=1}^\ell \|\xi^{m+\frac{1}{2}}\|_{L^3(\cT_h)}^3\notag\\
&-C\frac{k}{\epsilon^2}  (h^2+k+1)\sum_{m=1}^\ell\|\xi^{m}\|_{L^2(\cT_h)}^2-k^2\sum_{m=1}^\ell\|\xi^{m}\|_{L^2(\cT_h)}^2.\nonumber
\end{align}

Substitute the inequality above into \eqref{eq3.34}, and we get:
\begin{align}\label{eq3.50}
&\|\xi^\ell\|_{L^2(\cT_h)}^2 +\frac{2k}{\epsilon^2}\sum_{m=1}^\ell \bigl((\xi^{m})^2+(\xi^{m-1})^2, (\xi^{m}+\xi^{m-1})^2) \\
&+ 2k(1-\epsilon^2)\sum_{m=1}^\ell \left(  a_h(\xi^{m-\frac12},\xi^{m-\frac12})
+\frac{1}{\eps^2} \Bigl(f'\bigl(P_r^h u(t_{m-1})\bigr),(\xi^{m-\frac1 2})^2\Bigr)_{\cT_h}\right)\nonumber\\
 &+2k \sum_{m=1}^\ell \Bigl(f'\bigl(P_r^h u(t_{m-1})\bigr),(\xi^{m-\frac12})^2\Bigr)_{\cT_h}
\nonumber \\
&\leq \|\xi^{0}\|_{L^2(\cT_h)}^2 +Ch^{2\min\{r+1,s\}}\bigl(\, \|u\|_{H^1((0,T);H^s(\Ome))}^2 +\eps^{-4} \|u\|_{L^2((0,T);H^s(\Ome)}^2\bigr) \nonumber\\
&+\frac{C}{\epsilon^2}(h^8+k^4)+Ck^4[\epsilon^{-4}+\epsilon^{-2}+1]\nonumber\\
&+Ck(1+\frac{k\epsilon^2}{\epsilon^2}+\frac{h^2+k}{\epsilon^2})\sum_{m=1}^\ell \|\xi^{m}\|^2_{L^2(\cT_h)}+ C\frac{k}{\epsilon^2} \sum_{m=1}^\ell \|\xi^{m+\frac{1}{2}}\|_{L^3(\cT_h)}^3.\notag
\end{align}
\medskip
{\em Step 3}: In order to control the last two terms on the right-hand
side of \eqref{eq3.49}, we use the following Gagliardo-Nirenberg
inequality $\cite{23}$:
\[
\|v\|_{L^3(K)}^3\leq C\Bigl( \|\nab v\|_{L^2(K)}^{\frac{d}2}
\bigl\|v\bigr\|_{L^2(K)}^{\frac{6-d}2} +\|v\|_{L^2(K)}^3 \Bigr)
\qquad\forall K\in \cT_h,
\]
to get
\begin{align}\label{eq3.51}
\frac{Ck}{\eps^2} \sum_{m=1}^\ell \|\xi^m\|_{L^3(\cT_h)}^3
&\leq \eps^2\alpha k\sum_{m=1}^\ell \|\nab \xi^m\|_{L^2(\cT_h)}^2
+\eps^2 k\sum_{m=1}^\ell \|\xi^m\|_{L^2(\cT_h)}^2 \\
&\qquad
+C\eps^{-\frac{2(4+d)}{4-d}} k\sum_{m=1}^\ell\sum_{K\in \cT_h}
\bigl\|\xi^m\bigr\|_{L^2(K)}^{\frac{2(6-d)}{4-d}} \nonumber \\
&\leq \eps^2\alpha  k\sum_{m=1}^\ell \|\nab \xi^m\|_{L^2(\cT_h)}^2 \nonumber \\
&\qquad
+C\eps^{-\frac{2(4+d)}{4-d}} k\sum_{m=1}^\ell
\bigl\|\xi^m\bigr\|_{L^2(\cT_h)}^{\frac{2(6-d)}{4-d}}. \nonumber
\end{align}
Finally, for the third term on the left-hand side of the above inequality,
we utilize the discrete spectrum estimate \eqref{eq2.18} to
bound it from below as follows:
\begin{align}\label{eq3.52}
&2k(1-\epsilon^2)\sum_{m=1}^\ell \left(  a_h(\xi^{m-\frac12},\xi^{m-\frac12})
+\frac{1}{\eps^2} \Bigl(f'\bigl(P_r^h u(t_{m-1})\bigr),(\xi^{m-\frac1 2})^2\Bigr)_{\cT_h}\right) \\
&+4k \sum_{m=1}^\ell \Bigl(f'\bigl(P_r^h u(t_{m-1})\bigr),(\xi^{m-\frac12})^2\Bigr)_{\cT_h}\nonumber\\
&=2k(1-2\epsilon^2)\sum_{m=1}^\ell \left(  a_h(\xi^{m-\frac12},\xi^{m-\frac12})
+\frac{1}{\eps^2} \Bigl(f'\bigl(P_r^h u(t_{m-1})\bigr),(\xi^{m-\frac1 2})^2\Bigr)_{\cT_h}\right) \nonumber\\
&+2k\epsilon^2a_h(\xi^{m-\frac12},\xi^{m-\frac12})+ 4k \sum_{m=1}^\ell \Bigl(f'\bigl(P_r^h u(t_{m-1})\bigr),(\xi^{m-\frac12})^2\Bigr)_{\cT_h}\nonumber\\
&\geq -2(1-2\eps^2)c_0 k\sum_{m=1}^\ell \|\xi^m\|_{L^2(\cT_h)}^2 + 4\eps^2 \alpha k\sum_{m=1}^\ell \|\xi^{m-\frac{1}{2}}\|_{1,\DG}^2\notag\\
&-Ck \sum_{m=1}^\ell \|\xi^m\|_{L^2(\cT_h)}^2.\notag
\nonumber
\end{align}
\medskip
{\em Step 4}:
Substitute \eqref{eq3.51} and \eqref{eq3.52} into \eqref{eq3.50}, and we get the following:
\begin{align}\label{eq3.53}
&\|\xi^\ell\|_{L^2(\cT_h)}^2 + 3\eps^2 \alpha k\sum_{m=1}^\ell \|\xi^{m}\|_{1,\DG}^2+\frac{2k}{\epsilon^2}\sum_{m=1}^\ell \bigl((\xi^{m})^2+(\xi^{m-1})^2, (\xi^{m}+\xi^{m-1}))\\
&\leq Ck(1+\frac{k\epsilon^2}{\epsilon^2}+\frac{h^2+k}{\epsilon^2})\sum_{m=1}^\ell \|\xi^{m}\|^2_{L^2(\cT_h)}+C\eps^{-\frac{2(4+d)}{4-d}} k\sum_{m=1}^\ell
\bigl\|\xi^m\bigr\|_{L^2(\cT_h)}^{\frac{2(6-d)}{4-d}}\nonumber\\
&+\|\xi^{0}\|_{L^2(\cT_h)}^2 +Ch^{2\min\{r+1,s\}}\bigl(\, \|u\|_{H^1((0,T);H^s(\Ome))}^2 +\eps^{-4} \|u\|_{L^2((0,T);H^s(\Ome)}^2\bigr) \nonumber\\
&+\frac{C}{\epsilon^2}(h^8+k^4)+Ck^4[\epsilon^{-4}+\epsilon^{-2}+1].\nonumber\
\end{align}
Notice that on the right hand side, we need to choose the appropriate initial value $u^0_h$, so that $\|\xi^0\|_{L^2(\cT_h)}=O(h^{\min\{r+1,s\}})$
to maintain the optimal rate of convergence in $h$. Clearly,
both the $L^2$ and the elliptic projection of $u_0$ work.
and in the latter case, we get $\xi^0=0$.\\
It then follows from \eqref{eq2.7}, \eqref{eq2.9}, \eqref{eq2.12} and \eqref{eq3.53} that
\begin{align}\label{eq3.54}
&\|\xi^\ell\|_{L^2(\cT_h)}^2 + 3\eps^2 \alpha k\sum_{m=1}^\ell \|\xi^{m}\|_{1,\DG}^2+\frac{2k}{\epsilon^2}\sum_{m=1}^\ell \bigl((\xi^{m})^2+(\xi^{m-1})^2, (\xi^{m}+\xi^{m-1}))\\
&\leq Ck(1+\frac{k\epsilon^2}{\epsilon^2}+\frac{h^2+k}{\epsilon^2})\sum_{m=1}^\ell \|\xi^{m}\|^2_{L^2(\cT_h)}+C\eps^{-\frac{2(4+d)}{4-d}} k\sum_{m=1}^\ell
\bigl\|\xi^m\bigr\|_{L^2(\cT_h)}^{\frac{2(6-d)}{4-d}}\nonumber\\
& +C h^{2\min\{r+1,s\}}\eps^{-2(\sigma_1+2)}+\frac{C}{\epsilon^2}(h^8+k^4)+Ck^4[\epsilon^{-4}+\epsilon^{-2}+1].\nonumber\
\end{align}
Since $u^{\ell}$ can be written as
\begin{equation}\label{eq3.55}
u^{\ell}=k\mathop{\sum}\limits_{m=1}^\ell d_tu^m+u^0,
\end{equation}
then by \eqref{eq2.3} and \eqref{eq3.11}, we get
\begin{align}\label{eq3.56}
\|u^{\ell}\|_{L^2(\mathcal{T}_h)}
\leq k\mathop{\sum}\limits_{m=1}^\ell \|d_tu^m\|_{L^2(\mathcal{T}_h)}
+\|u^0\|_{L^2(\mathcal{T}_h)}
\leq C\eps^{-2\sigma_1}.
\end{align}
By the boundedness of the projection, we have
\begin{equation}\label{eq3.57}
\|\xi^\ell\|_{L^2(\mathcal{T}_h)}^2\leq C\eps^{-2\sigma_1}.
\end{equation}
Then the above inequality is equivalent to the form below:
\begin{equation}\label{eq3.58}
\|\xi^\ell\|_{L^2(\cT_h)}^2 + k\sum_{m=1}^\ell
3\eps^2 \alpha \|\xi^m\|_{1,\DG}^2 \leq H_1+H_2,
\end{equation}
where
\begin{align}\label{eq3.59}
H_1:&=Ck(1+\frac{k\epsilon^2}{\epsilon^2}+\frac{h^2+k}{\epsilon^2})\sum_{m=1}^{\ell-1} \|\xi^{m}\|^2_{L^2(\cT_h)}+C\eps^{-\frac{2(4+d)}{4-d}} k\sum_{m=1}^{\ell-1}
\bigl\|\xi^m\bigr\|_{L^2(\cT_h)}^{\frac{2(6-d)}{4-d}}\\ \nonumber
&\qquad\qquad
+C h^{2\min\{r+1,s\}}\eps^{-2(\sigma_1+2)}+\frac{C}{\epsilon^2}(h^8+k^4)+Ck^4[\epsilon^{-4}+\epsilon^{-2}+1],\nonumber
\end{align}
\begin{align}\label{eq3.60}
H_2:&=Ck(1+\frac{k\epsilon^2}{\epsilon^2}+\frac{h^2+k}{\epsilon^2})\|\xi^{\ell}\|_{L^2(\cT_h)}^2 + C\eps^{-\frac{2(4+d)}{4-d}} k
\bigl\|\xi^{\ell}\bigr\|_{L^2(\cT_h)}^{\frac{2(6-d)}{4-d}}.
\end{align}
It is easy to check that
\begin{equation}\label{eq3.61}
H_2<\frac12\|\xi^\ell\|_{L^2(\mathcal{T}_h)}^2 ,
\qquad\mbox{provided that}\quad k<A(\epsilon).
\end{equation}
By \eqref{eq3.58} we have
\begin{align}\label{eq3.62}
& \|\xi^\ell\|_{L^2(\cT_h)}^2 + k\sum_{m=1}^\ell
3\eps^2 \alpha \|\xi^m\|_{1,\DG}^2 \leq2H_1\\\nonumber
&
\leq 2Ck(1+\frac{k\epsilon^2}{\epsilon^2}+2\frac{h^2+k}{\epsilon^2})\sum_{m=1}^{\ell-1} \|\xi^{m}\|^2_{L^2(\cT_h)}+2C\eps^{-\frac{2(4+d)}{4-d}} k\sum_{m=1}^{\ell-1}
\bigl\|\xi^m\bigr\|_{L^2(\cT_h)}^{\frac{2(6-d)}{4-d}}\nonumber\\
&
+2C h^{2\min\{r+1,s\}}\eps^{-2(\sigma_1+2)}+2\frac{C}{\epsilon^2}(h^8+k^4)+2Ck^4[\epsilon^{-4}+\epsilon^{-2}+1]\nonumber\\
&
\leq Ck(1+\frac{k\epsilon^2}{\epsilon^2}+\frac{h^2+k}{\epsilon^2})\sum_{m=1}^{\ell-1} \|\xi^{m}\|^2_{L^2(\cT_h)}+C\eps^{-\frac{2(4+d)}{4-d}} k\sum_{m=1}^{\ell-1}
\bigl\|\xi^m\bigr\|_{L^2(\cT_h)}^{\frac{2(6-d)}{4-d}}\nonumber\\
&
+C h^{2\min\{r+1,s\}}\eps^{-2(\sigma_1+2)}+\frac{C}{\epsilon^2}(h^8+k^4)+Ck^4[\epsilon^{-4}+\epsilon^{-2}+1].\nonumber
\end{align}
Let $d_{\ell}\geq 0$ be the slack variable such that
\begin{align}\label{eq3.63}
& \|\xi^\ell\|_{L^2(\mathcal{T}_h)}^2 + k\sum_{m=1}^\ell
3\eps^2 \alpha \|\xi^m\|_{1,\DG}^2 +d_{\ell} \\
&
=Ck(1+\frac{k\epsilon^2}{\epsilon^2}+\frac{h^2+k}{\epsilon^2})\sum_{m=1}^{\ell-1} \|\xi^{m}\|^2_{L^2(\cT_h)}+C\eps^{-\frac{2(4+d)}{4-d}} k\sum_{m=1}^{\ell-1}
\bigl\|\xi^m\bigr\|_{L^2(\cT_h)}^{\frac{2(6-d)}{4-d}} \nonumber\\
&
+C h^{2\min\{r+1,s\}}\eps^{-2(\sigma_1+2)}+\frac{C}{\epsilon^2}(h^8+k^4)+Ck^4[\epsilon^{-4}+\epsilon^{-2}+1].\nonumber
\end{align}
and define for $\ell\geq1$
\begin{align}\label{eq3.64}
S_{\ell+1}:&=  \|\xi^\ell\|_{L^2(\mathcal{T}_h)}^2 + k\sum_{m=1}^\ell
 3\eps^2 \alpha \|\xi^m\|_{1,\DG}^2+d_{\ell},
\end{align}
\begin{align}\label{eq3.65}
S_{1}:&=C h^{2\min\{r+1,s\}}\eps^{-2(\sigma_1+2)}+\frac{C}{\epsilon^2}(h^8+k^4)+Ck^4[\epsilon^{-4}+\epsilon^{-2}+1],
\end{align}
then we have
\begin{equation}\label{eq3.66}
S_{\ell+1}-S_{\ell}\leq C(1+\frac{k\epsilon^2}{\epsilon^2}+\frac{h^2+k}{\epsilon^2}) kS_{\ell}+C\eps^{-\frac{2(4+d)}{4-d}} kS_{\ell}^{\frac{6-d}{4-d}}\qquad\text{for}\ \ell\geq1.
\end{equation}
Applying Lemma \ref{lem2.4} to $\{S_\ell\}_{\ell\geq 1}$ defined above,
we obtain for $\ell\geq1$
\begin{equation}\label{eq3.67}
S_{\ell}\leq a^{-1}_{\ell}\Bigg\{S^{-\frac{2}{4-d}}_{1}-\frac{2Ck}{4-d}
\sum_{s=1}^{\ell-1}\eps^{-\frac{2(4+d)}{4-d}} a^{-\frac{2}{4-d}}_{s+1}\Bigg\}^{-\frac{4-d}{2}}.
\end{equation}
provided that
\begin{equation}\label{eq3.68}
\frac12 S^{-\frac{2}{4-d}}_{1}-\frac{2Ck}{4-d}\sum_{s=1}^{\ell-1} \eps^{-\frac{2(4+d)}{4-d}}
a^{-\frac{2}{4-d}}_{s+1}>0.
\end{equation}
We note that $a_s\, (1\leq s\leq \ell)$ are all bounded as $k\rightarrow0$,
therefore, \eqref{eq3.68} holds under the mesh constraint stated in the theorem.
It follows from \eqref{eq3.66} and \eqref{eq3.67} that
\begin{equation}\label{eq3.69}
S_{\ell}\leq 2a_\ell^{-1} S_1
\leq Ck^4\eps^{-2(\sigma_1+2)}+C h^{2\min\{r+1,s\}}\eps^{-2(\sigma_1+2)}.
\end{equation}

Finally, using the above estimate and the properties of the operator $P^h_r$
we obtain \eqref{eq3.18} and \eqref{eq3.19}. The estimate \eqref{eq3.20} follows
from \eqref{eq3.19} and the inverse inequality bounding the $L^\infty$-norm by the
$L^2$-norm and \eqref{eq2.21}.  The proof is complete.

\section{Convergence of the numerical interface to the mean curvature flow} \label{sec-4}
In this section, we prove the rate of convergence of the numerical interface to its limit geometric interface of the Allen-Cahn  equation. This convergence theory is based on the maximum norm error estimates, which is proven above. The rate of convergence can be proven by the sharper error estimates, which is the negative polynomial function of the interaction length $\epsilon$ \cite{15, 17,3}. It can't be proven if the coarse error estimate, which is the exponential function of $\epsilon$, is used.

For all the DG problem, the the zero-level set of $u_h^n$ may not be well defined since the zero-level set may not be continuous. Therefore, we introduce the finite element approximation
$\widehat{u}_h^m$ of the DG solution $u_h^m$ It is defined by using the averaged degrees of freedom
of $u_h^n$ as the degrees of freedom for determining $\widehat{u}_h^m$ (cf. $\cite{24}$). We get the following results $\cite{24}$.

\begin{theorem}\label{lem4.1}
Let $\mathcal{T}_h$ be a conforming mesh consisting of
triangles when $d=2$, and tetrahedra when $d=3$. For $v_h\in V_h$, let
$\widehat{v}_h$ be the finite element approximation  of $v_h$ as
defined above. Then for any $v_h\in V_h$ and $i=0,1$ there holds
\begin{align}\label{eqn_KP}
\sum_{K\in\mathcal{T}_h} \|v_h-\widehat{v}_h\|_{H^i(K)}^2
\leq C \sum_{e\in\cE_h^I} h^{1-2i}_e \|[v_h]\|_{L^2(e)}^2,
\end{align}
where $C>0$ is a constant independent of $h$ and $v_h$ but may depend on $r$ and the minimal
angle $\theta_0$ of the triangles in $\mathcal{T}_h$.
\end{theorem}
\\
Using the above approximation result we can show that the error estimates
of Theorem \ref{thm3.5} also hold for $\widehat{u}_h^n$.
\begin{theorem}\label{lem4.2}
Let $u_h^{m}$ denote the solution of the DG scheme \eqref{eq3.1}--\eqref{eq3.4}
and $\widehat{u}_h^{m}$ denote its finite element approximation as defined above. Then
under the assumptions of Theorem \ref{thm3.5} the error estimates for $u_h^m$ given in
Theorem \ref{thm3.5} are still valid for $\widehat{u}_h^{m}$, in particular, there holds
\begin{align}\label{eq3.36bx}
\max_{0\leq m\leq M} \|u(t_m)-\widehat{u}_h^m\|_{L^\infty(\cT_h)}
&\leq C h^{\min\{r+1,s\}} |\ln h|^{\overline{r}} \eps^{-\gamma}\\
&\qquad
+Ch^{-\frac{d}{2}}(k^2+h^{\min\{r+1,s\}})\eps^{-(\sigma_1+2)}.  \nonumber
\end{align}

\end{theorem}


Proof: We only give a proof for \eqref{eq3.36bx} because other estimates can
be proved likewise.  By the triangle inequality we have
\begin{align}\label{eq3.36by}
\|u(t_m)-\widehat{u}_h^m\|_{L^\infty(\cT_h)}
\leq \|u(t_m)-u_h^m\|_{L^\infty(\cT_h)}+ \|u_h^m-\widehat{u}_h^m\|_{L^\infty(\cT_h)}.
\end{align}
Hence, it suffices to show that the second term on the right-hand side
is an equal or higher order term compared to the first one.

Let $u^I(t)$ denote the finite element interpolation of $u(t)$ into $S_h$.
It follows from \eqref{eqn_KP} and the trace inequality that
\begin{align}\label{eq3.36bz}
\|u_h^m-\widehat{u}_h^m\|_{L^2(\cT_h)}^2
&\leq C\sum_{e\in \cE_h^I} h_e \|[u_h^m]\|_{L^2(e)}^2 \\
&= C\sum_{e\in \cE_h^I} h_e \|[u_h^m-u^I(t_m)]\|_{L^2(e)}^2 \nonumber\\
& \leq C\sum_{K\in \cT_h} h_e h_K^{-1}\|u_h^m-u^I(t_m)\|_{L^2(K)}^2 \nonumber \\
&\leq C \bigl( \|u_h^m-u(t_m)\|_{L^2(\cT_h)}^2 + \|u(t_m)- u^I(t_m)\|_{L^2(\cT_h)}^2 \bigr). \nonumber
\end{align}
Substituting \eqref{eq3.36bz} into \eqref{eq3.36by} after using the inverse inequality yields
\begin{align*}
&\|u(t_m)-\widehat{u}_h^m\|_{L^\infty(\cT_h)}
\leq \|u(t_m)-u_h^m\|_{L^\infty(\cT_h)}+ C h^{-\frac{d}2} \|u_h^m-\widehat{u}_h^m\|_{L^2(\cT_h)}\\
&\qquad\quad
\leq \|u(t_m)-u_h^m\|_{L^\infty(\cT_h)}  \nonumber \\
&\qquad\qquad
+ Ch^{-\frac{d}2}  \bigl( \|u_h^m-u(t_m)\|_{L^2(\cT_h)} + \|u(t_m)- u^I(t_m)\|_{L^2(\cT_h)} \bigr),\nonumber
\end{align*}
which together with \eqref{eq3.18} implies the desired estimate \eqref{eq3.36bx}. The proof is complete.

We are now ready to state the main theorem of this section.

\begin{theorem}\label{thm4.3}
Let $\{\Gamma_t\}$ denote the (generalized) mean curvature flow
defined in $\cite{25}$, that is, $\Gamma_t$ is the zero-level set of
the solution $w$ of the following initial value problem:
\begin{alignat}{2}\label{eq4.1}
w_t &=\Delta w-\frac{D^2wDw\cdot Dw}{|Dw|^2} &&\qquad\mbox{in }
\mathbf{R}^d\times (0,\infty), \\
w(\cdot,0) &=w_0(\cdot) &&\qquad\mbox{in } \mathbf{R}^d. \label{eq4.2}
\end{alignat}
Let $u^{\epsilon,h,k}$ denote the piecewise linear interpolation in time
of the numerical solution $\{\widehat{u}_h^m\}$ defined by
\begin{equation}\label{eq4.3}
u^{\epsilon,h,k}(x,t):=\frac{t-t_m}{k}\widehat{u}_h^{m+1}(x)+\frac{t_{m+1}-t}{k}
\widehat{u}_h^{m}(x), \quad t_m\leq t\leq t_{m+1}
\end{equation}
for $0\leq m\leq M-1$. Let $\{\Gamma_t^{\epsilon,h,k}\}$ denote the zero-level
set of $u^{\epsilon,h,k}$, namely,
\begin{equation}\label{eq4.4}
\Gamma_t^{\epsilon,h,k}=\{x\in \Omega;\, u^{\epsilon,h,k}(x,t)=0\}.
\end{equation}
Suppose $\Gamma_0=\{x\in \overline{\Omega};u_0(x)=0\}$ is a smooth hypersurface
compactly contained in $\Omega$, and $k=O(h^2)$. Let $t_*$ be
the first time at which the mean curvature flow develops a singularity, then
there exists a constant $\epsilon_1>0$ such that for all
$\epsilon\in(0,\epsilon_1)$ and $ 0<t<t_*$ there holds
\[
\sup_{x\in\Gamma_t^{\epsilon,h,k}}\{\mbox{\rm dist}(x,\Gamma_t)\}
\leq C\epsilon^2|\ln\,\epsilon|^2.
\]

\end{theorem}
Proof: We note that since $u^{\epsilon,h,k}(x,t)$ is continuous in both $t$ and $x$, then
$\Gamma_t^{\epsilon,h,k}$ is well defined.
Let $I_t$ and $O_t$ denote the inside and the outside of $\Gamma_t$ defined by
\begin{equation}\label{eq4.5}
I_t:=\{ x\in \mathbf{R}^d;\, w(x,t)>0\}, \qquad O_t:=\{ x\in \mathbf{R}^d;\, w(x,t)<0\}.
\end{equation}

Let $d(x,t)$ denote the signed distance function to $\Gamma_t$ which is positive
in $I_t$ and negative in $O_t$. By Theorem 6.1 of $\cite{26}$, there exist
$\widehat{\epsilon}_1>0$ and $\widehat{C}_1>0$ such that for all $t\geq  0$
and $\epsilon\in(0,\widehat{\epsilon}_1)$ there hold
\begin{alignat}{2}\label{eq4.6}
u_{\epsilon}(x,t) &\geq 1-\epsilon
&&\qquad\forall x\in\{x\in\overline{\Omega};\, d(x,t)\geq \widehat{C}_1\epsilon^2|
\ln\,\epsilon|^2\},\\
u_{\epsilon}(x,t) &\leq -1+\epsilon
&&\qquad\forall x\in\{x\in\overline{\Omega};\, d(x,t)\leq -\widehat{C}_1\epsilon^2|
\ln\,\epsilon|^2\}. \label{eq4.7}
\end{alignat}

Since for any fixed $x\in\Gamma_t^{\epsilon,h,k}$, $u^{\epsilon,h,k}(x,t)=0$,
by \eqref{eq3.36bx} with $k=O(h^2)$, we have
\begin{align*}
|u^{\epsilon}(x,t)| &=|u^{\epsilon}(x,t)-u^{\epsilon,h,k}(x,t)| \\
&\leq \tilde{C} \Bigl( h^{\min\{r+1,s\}} |\ln h|^{\overline{r}} \eps^{-\gamma}
+h^{-\frac{d}{2}}(k+h^{\min\{r+1,s\}})\eps^{-(\sigma_1+2)} \Bigr).
\end{align*}
Then there exists $\widetilde{\epsilon}_1>0$ such that for
$\epsilon\in(0,\widetilde{\epsilon}_1)$ there holds
\begin{equation}\label{eq4.8}
|u^{\epsilon}(x,t)|<1-\epsilon.
\end{equation}
Therefore, the assertion follows from setting
$\epsilon_1=\min\{\widehat{\epsilon}_1,\widetilde{\epsilon}_1\}$.
The proof is complete.
\section{Numerical experiments}\label{sec-5}

In this section, we provide two two-dimensional numerical experiments to gauge the accuracy and reliability of the
 fully discrete IPDG method developed in the previous sections. We use a square domain $\Omega=[-1,1]\times[-1,1] \subset\mathbf{R}^2$, and $u_0(x)=\mbox{tanh}(\frac{d_0(x)}{\sqrt{2}\epsilon})$,
where $d_0(x)$ stands for the signed distance from $x$ to the initial curve $\Gamma_0$ See the details for similar numerical setting in \cite{feng2014multiphysics, Feng_Li15, feng2015analysis,16,li2015numerical, xu2016convex}.

The first test uses the smooth initial curves $\Gamma_0$, hence the requirements for $u_0$ are satisfied. Consequently, the results established in this paper apply to the test example.
In the test we first verify the spatial rate of convergence given in \eqref{eq3.18} and \eqref{eq3.20}. We then compute
the evolution of the zero-level set of the solution of the Allen-Cahn problem with
$\epsilon= 0.025$ and at various time instances.

{\bf Test 1} Consider the Allen-Cahn problem with the following initial condition:
$$
u_0(x)=\left\{\begin{array}{ll}\mbox{tanh}(\frac{d(x)}{\sqrt{2}\epsilon}), & \mbox{if}\ \frac{x_1^2}{0.36}+\frac{x_2^2}{0.04}\geq1,\\
                              \mbox{tanh}(\frac{-d(x)}{\sqrt{2}\epsilon}),& \mbox{if}\ \frac{x_1^2}{0.36}+\frac{x_2^2}{0.04}<1,
 \end{array}
 \right.
$$
here $d(x)$ stands for the distance function to the ellipse $\frac{x_1^2}{0.36}+\frac{x_2^2}{0.04}=1$.

\begin{center}{{ Table 5.1.\ \ Spatial errors and convergence rates\\}}
\vspace{2mm} {
\begin{small}
\begin{tabular}{ccccc}
\hline
\hspace{0.3cm}$h$\hspace{0.3cm}&$L^\infty(L^2)$ error \hspace{0.3cm}&$L^\infty(L^2)$ order\hspace{0.3cm}&$L^2(H^1)$ error\hspace{0.3cm}& $L^2(H^1)$ order\hspace{0.3cm}
\\
\hline
\hspace{0.3cm}$\sqrt{2}/10$\hspace{0.3cm}&0.02451 \hspace{0.3cm}&\hspace{0.3cm}&0.34216\hspace{0.3cm}& \hspace{0.3cm}\\
\hspace{0.3cm}$\sqrt{2}/20$\hspace{0.3cm}&0.00539 \hspace{0.3cm}&2.1850\hspace{0.3cm}&0.17258\hspace{0.3cm}& 0.9874\hspace{0.3cm}\\
\hspace{0.3cm}$\sqrt{2}/40$\hspace{0.3cm}&0.00142 \hspace{0.3cm}&1.9244\hspace{0.3cm}&0.08394\hspace{0.3cm}& 1.0398\hspace{0.3cm}\\
\hspace{0.3cm}$\sqrt{2}/80$\hspace{0.3cm}&0.00036 \hspace{0.3cm}&1.9798\hspace{0.3cm}&0.04172\hspace{0.3cm}&1.0086\hspace{0.3cm}\\
\hline
\end{tabular}\end{small}
}
\end{center}
Table 5.1 shows the spatial $L^2$ and $H^1$-norm errors and convergence rates,
which are consistent with what are proved for the linear element in the convergence theorem.
\begin{figure}[th]
\centering
\includegraphics[height=1.8in,width=2.4in]{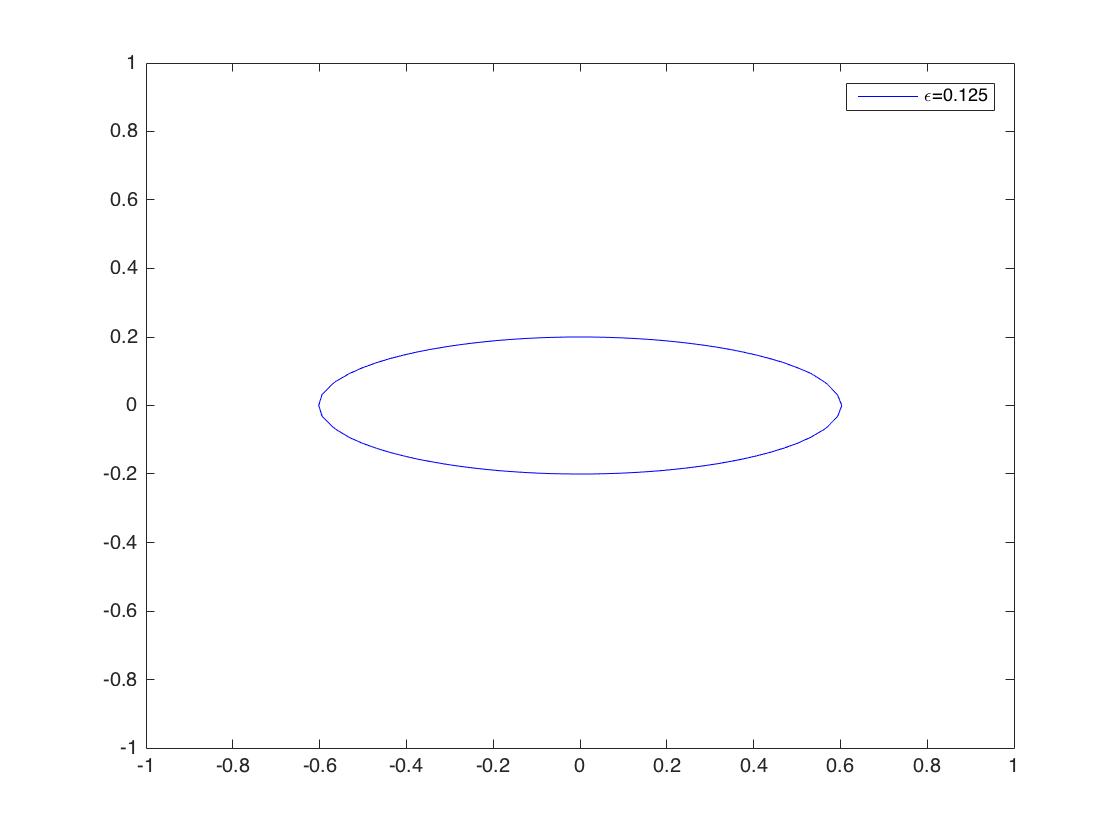}
\includegraphics[height=1.8in,width=2.4in]{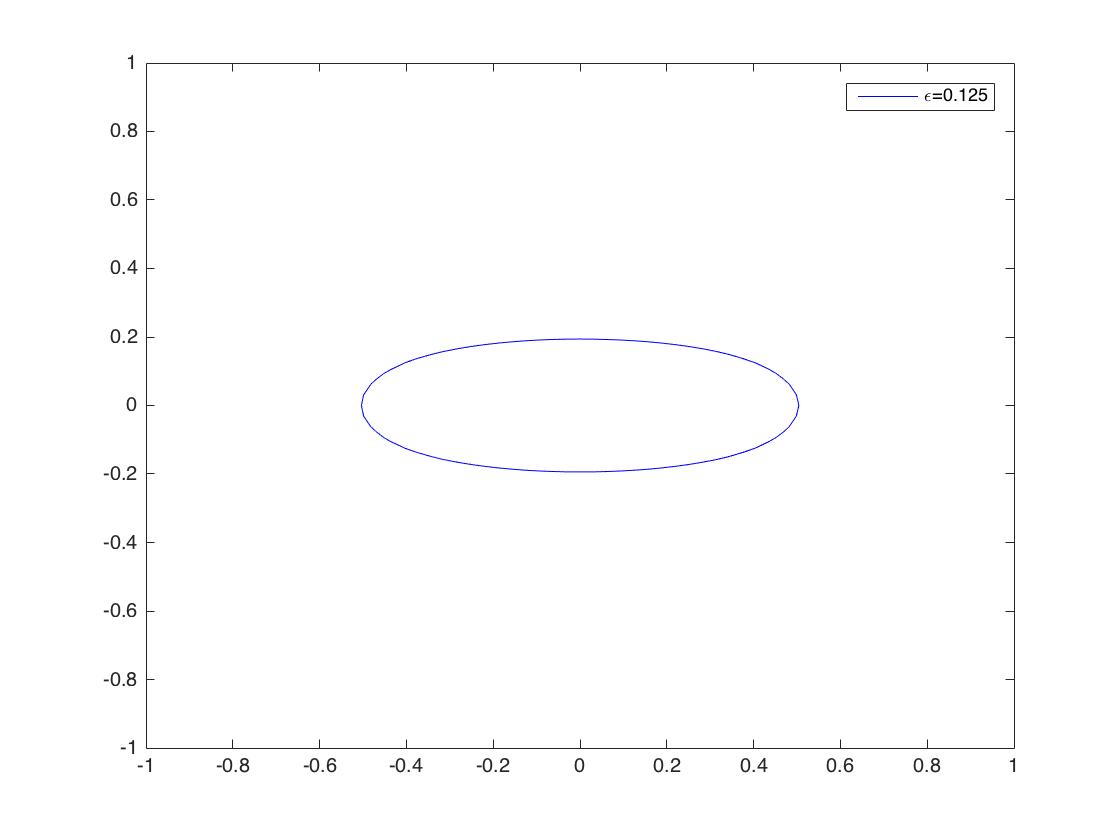}

\includegraphics[height=1.8in,width=2.4in]{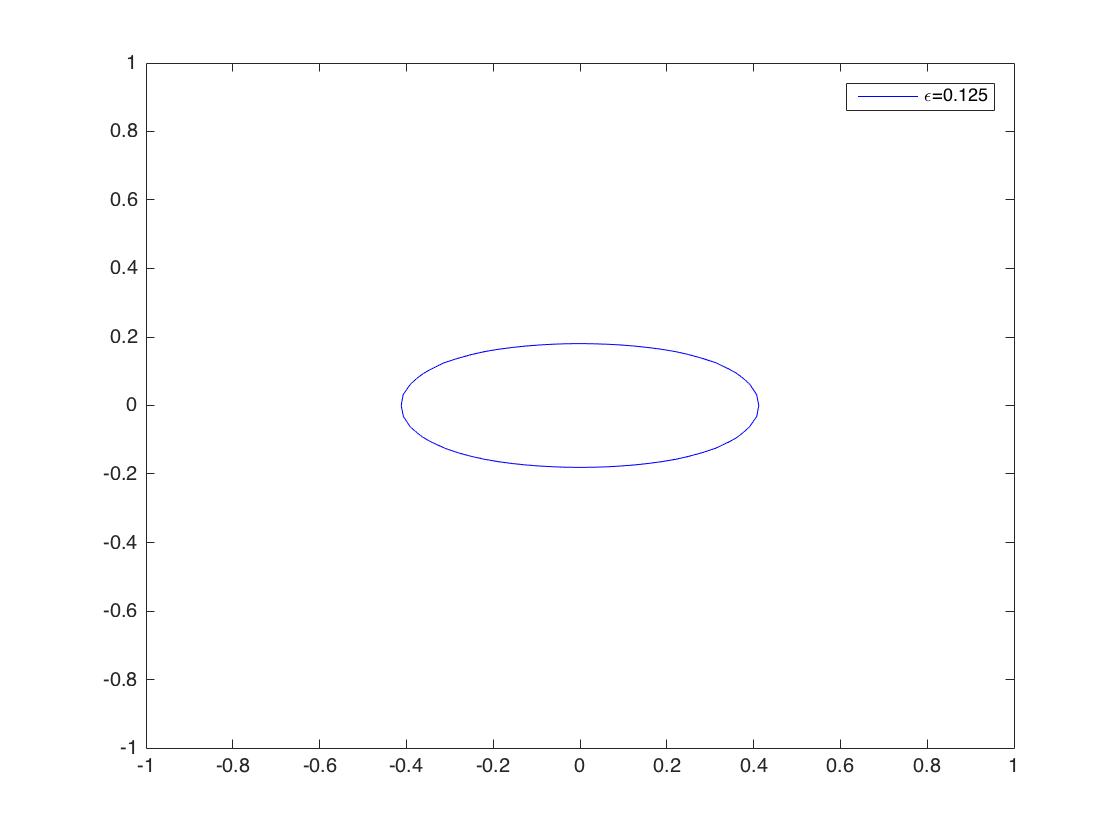}
\includegraphics[height=1.8in,width=2.4in]{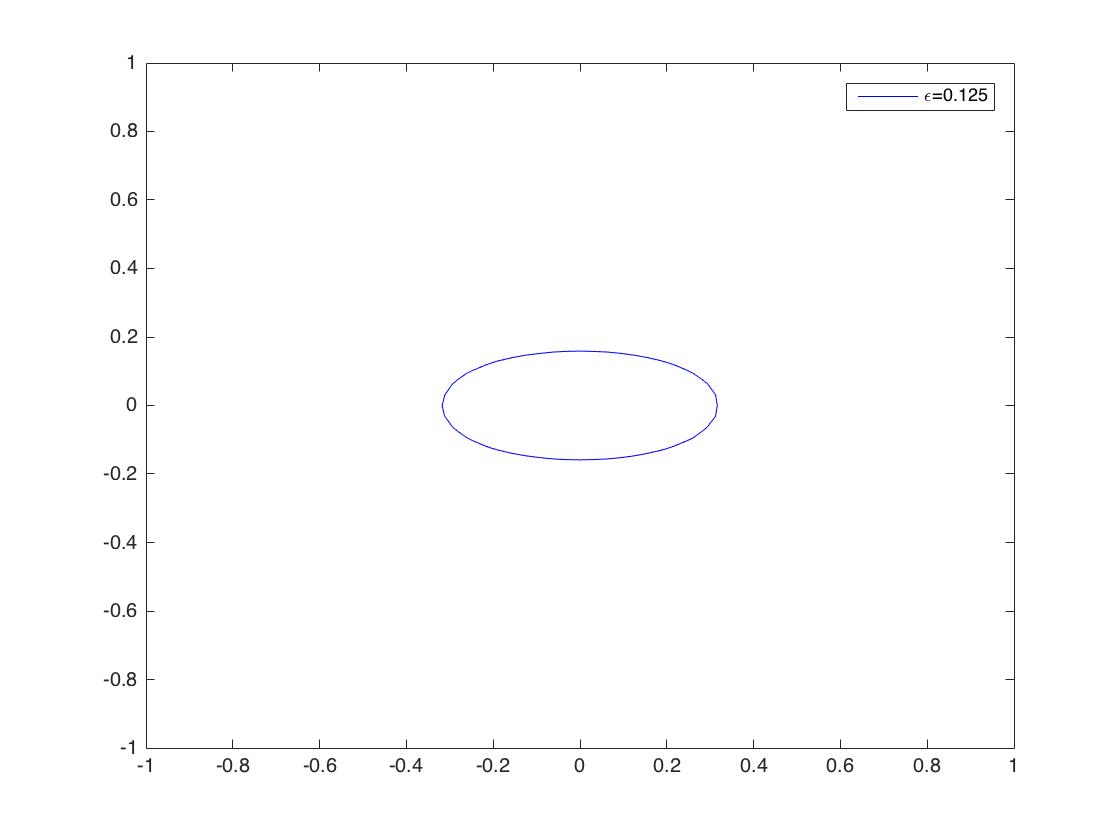}

\includegraphics[height=1.8in,width=2.4in]{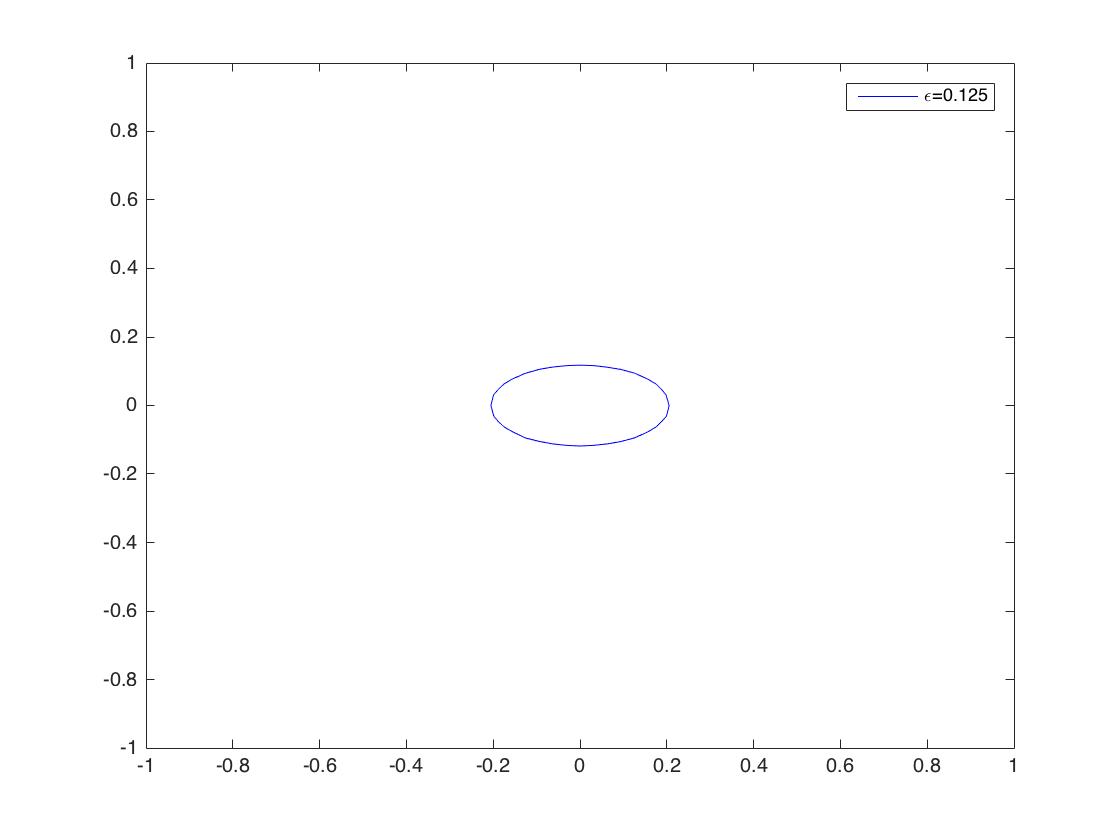}
\includegraphics[height=1.8in,width=2.4in]{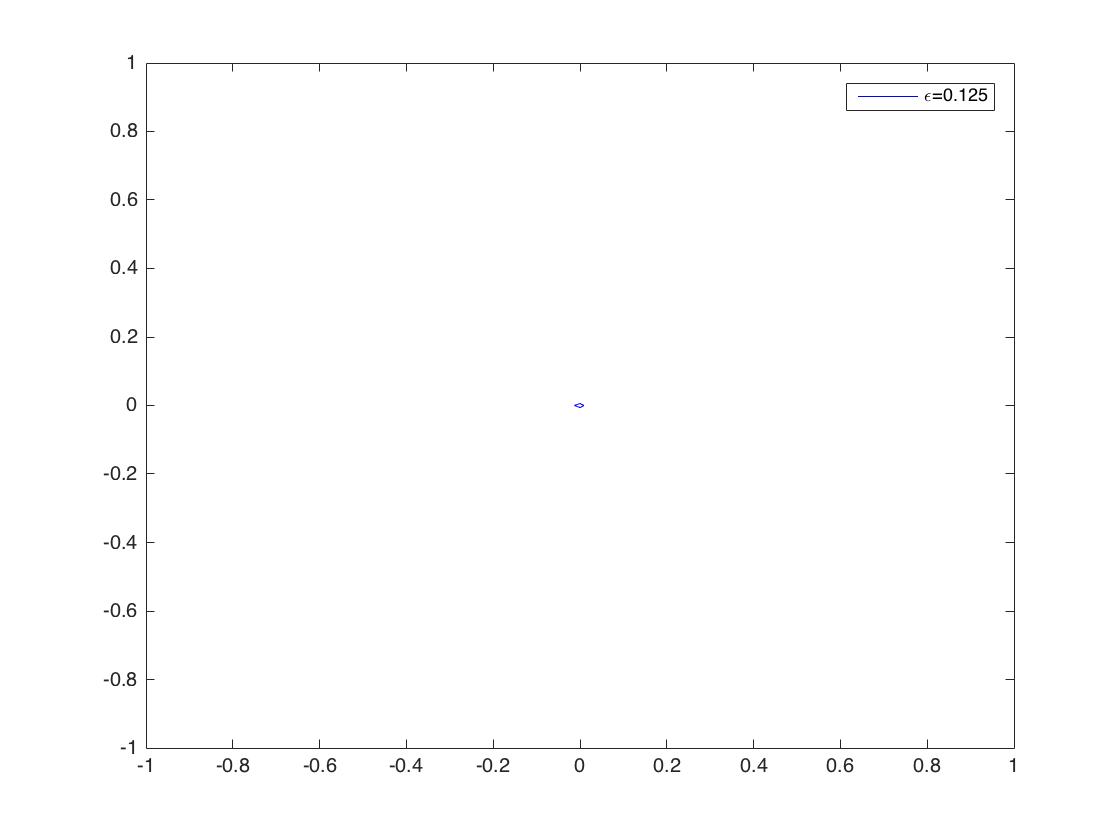}

\caption{Test 1: Snapshots of the zero-level set of $u^{\epsilon,h,k}$ at time $t=0,
0.8\times 10^{-2}, 1.6\times 10^{-2}, 2.4\times 10^{-2},3.2\times 10^{-2},.3.8\times 10^{-2}$ and $\epsilon=0.125$.}
\end{figure}

Figure 5.1 displays six snapshots of the zero-level set of the numerical solution $u^{\epsilon,h,k}$ with $\epsilon=0.125$.
We observe that as $\epsilon$ is small enough the zero-level set converges to the mean curvature flow $\Gamma_t$
as time goes on.

The second test has non-smooth curve with $u_0$ defined below.This initial condition does not satisfy the assumptions in the spetrum estimate, but we can still numerically validate the convergence of the solution to the mean curvature flow.\\
{\bf Test 2} Consider the Allen-Cahn problem with the following initial condition:
\begin{equation*}
 u_0(x)=\begin{cases}
 \tanh(\frac{1}{\sqrt{2}\eps}(\text{min}\{d_1(x),d_2(x)\})), & \text{if}\
\frac{x_1^2}{0.36}+\frac{x_2^2}{0.04}\geq1,\frac{x_1^2}{0.04}+\frac{x_2^2}{0.36}\geq1,\\
 &\ \text{or}\ \frac{x_1^2}{0.36}+\frac{x_2^2}{0.04}\leq1,\frac{x_1^2}{0.04}+\frac{x_2^2}{0.04}\leq1,\\
  \tanh(\frac{-1}{\sqrt{2}\eps}(\text{min}\{d_1(x),d_2(x)\})), & \text{if}\
\frac{x_1^2}{0.36}+\frac{x_2^2}{0.04}<1,\frac{x_1^2}{0.04}+\frac{x_2^2}{0.36}>1,\\
 &\ \text{or}\ \frac{x_1^2}{0.36}+\frac{x_2^2}{0.04}>1,\frac{x_1^2}{0.04}+\frac{x_2^2}{0.36}<1.\\
\end{cases}
\end{equation*}

here $d_1(x)$  and $d_2(x)$stands for the distance function to the ellipses $\frac{x_1^2}{0.36}+\frac{x_2^2}{0.04}=1$ and $\frac{x_1^2}{0.04}+\frac{x_2^2}{0.36}=1$ respectively.

\begin{center}{{ Table 5.2.\ \ Spatial errors and convergence rates\\}}
\vspace{2mm} {
\begin{small}
\begin{tabular}{ccccc}
\hline
\hspace{0.3cm}$h$\hspace{0.3cm}&$L^\infty(L^2)$ error \hspace{0.3cm}&$L^\infty(L^2)$ order\hspace{0.3cm}&$L^2(H^1)$ error\hspace{0.3cm}& $L^2(H^1)$ order\hspace{0.3cm}
\\
\hline
\hspace{0.3cm}$\sqrt{2}/10$\hspace{0.3cm}&0.01032 \hspace{0.3cm}&\hspace{0.3cm}&0.08325\hspace{0.3cm}& \hspace{0.3cm}\\
\hspace{0.3cm}$\sqrt{2}/20$\hspace{0.3cm}&0.00256 \hspace{0.3cm}&2.0098\hspace{0.3cm}&0.03851\hspace{0.3cm}& 1.1123\hspace{0.3cm}\\
\hspace{0.3cm}$\sqrt{2}/40$\hspace{0.3cm}&0.00075 \hspace{0.3cm}&1.7638\hspace{0.3cm}&0.01888\hspace{0.3cm}& 1.0283\hspace{0.3cm}\\
\hspace{0.3cm}$\sqrt{2}/80$\hspace{0.3cm}&0.00022 \hspace{0.3cm}&1.9836\hspace{0.3cm}&0.00939\hspace{0.3cm}&1.0069\hspace{0.3cm}\\
\hline
\end{tabular}\end{small}
}
\end{center}
Table 5.2 shows the spatial $L^2$ and $H^1$-norm errors and convergence rates,
which are consistent with what are proved for the linear element in the convergence theorem.
\begin{figure}[ht]
\centering
\includegraphics[height=1.8in,width=2.4in]{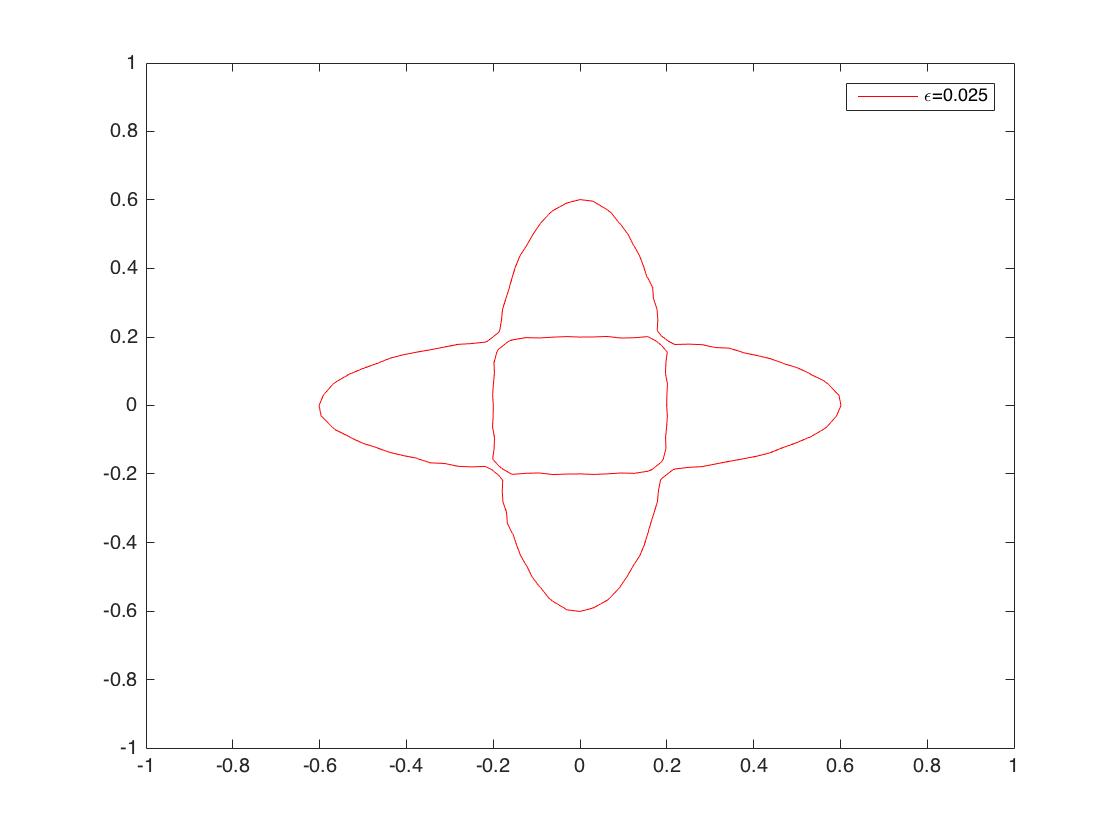}
\includegraphics[height=1.8in,width=2.4in]{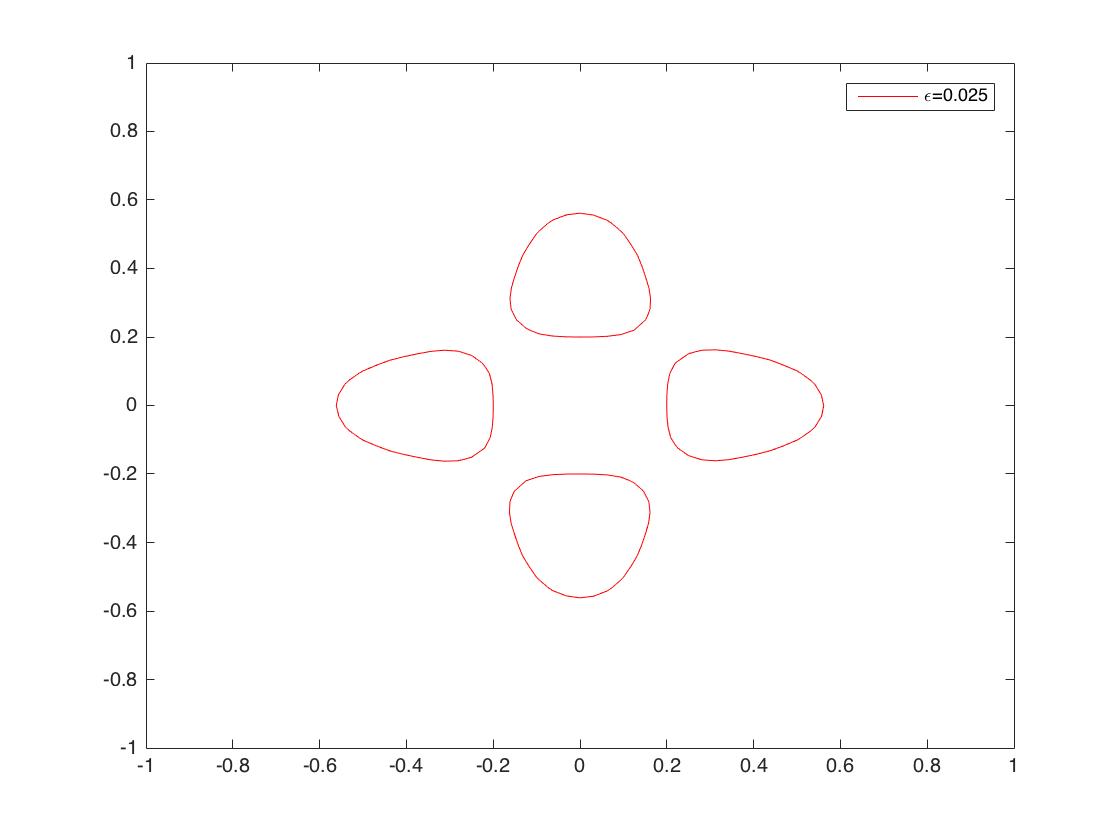}

\includegraphics[height=1.8in,width=2.4in]{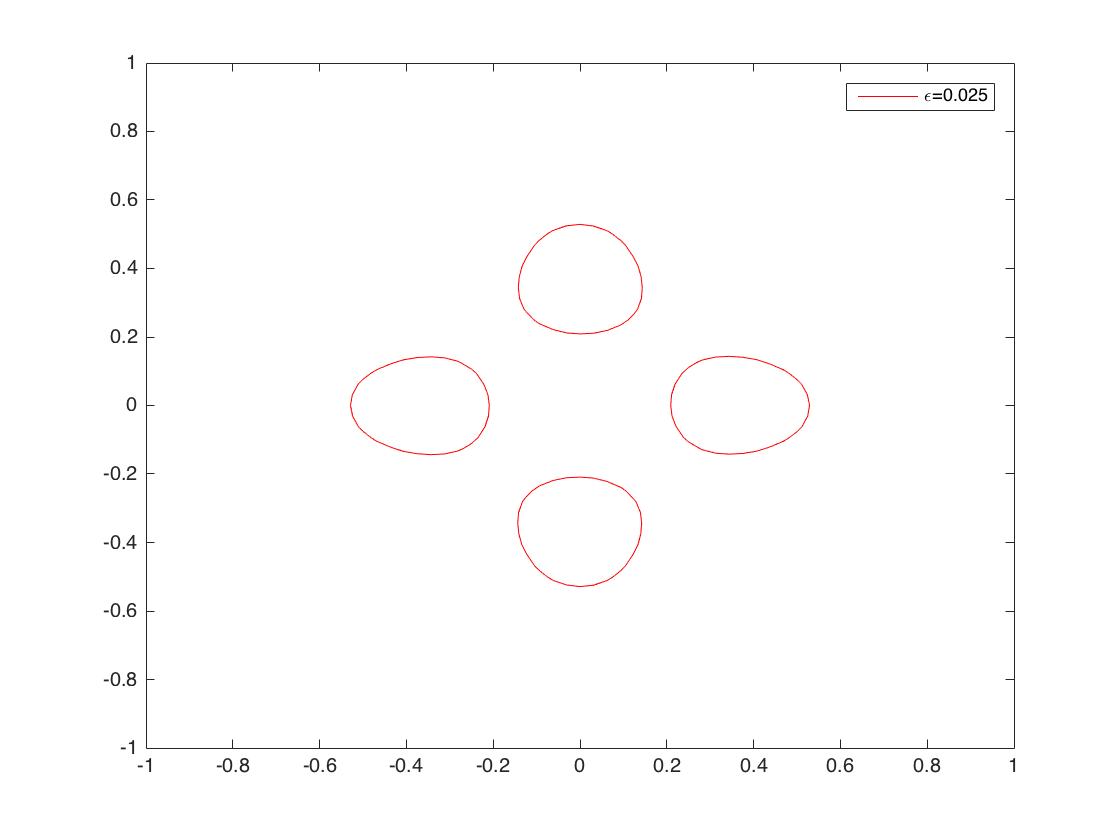}
\includegraphics[height=1.8in,width=2.4in]{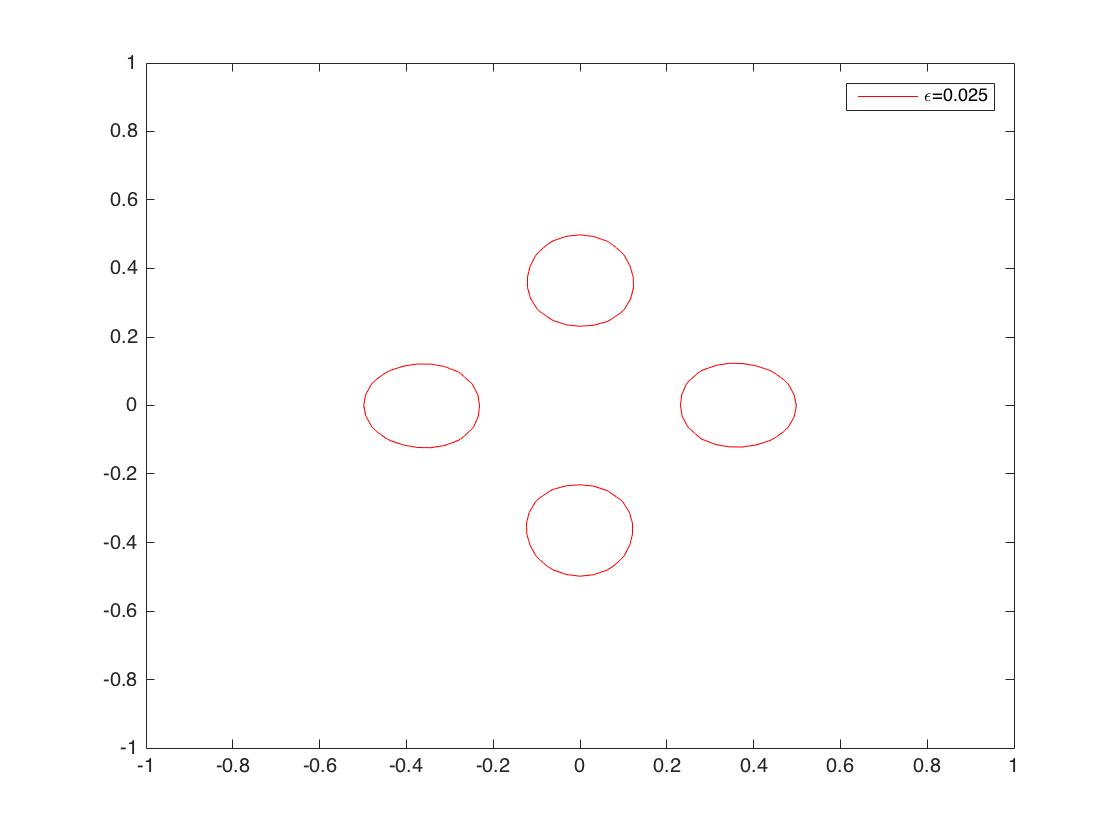}

\includegraphics[height=1.8in,width=2.4in]{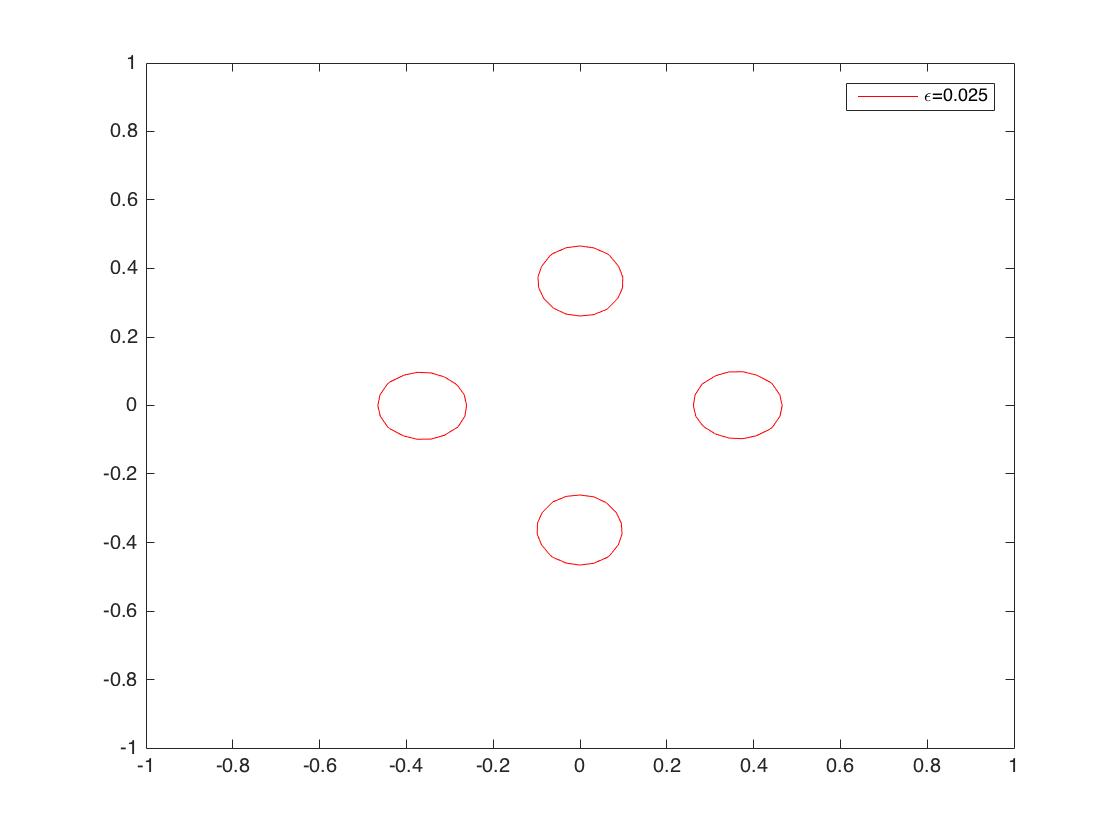}
\includegraphics[height=1.8in,width=2.4in]{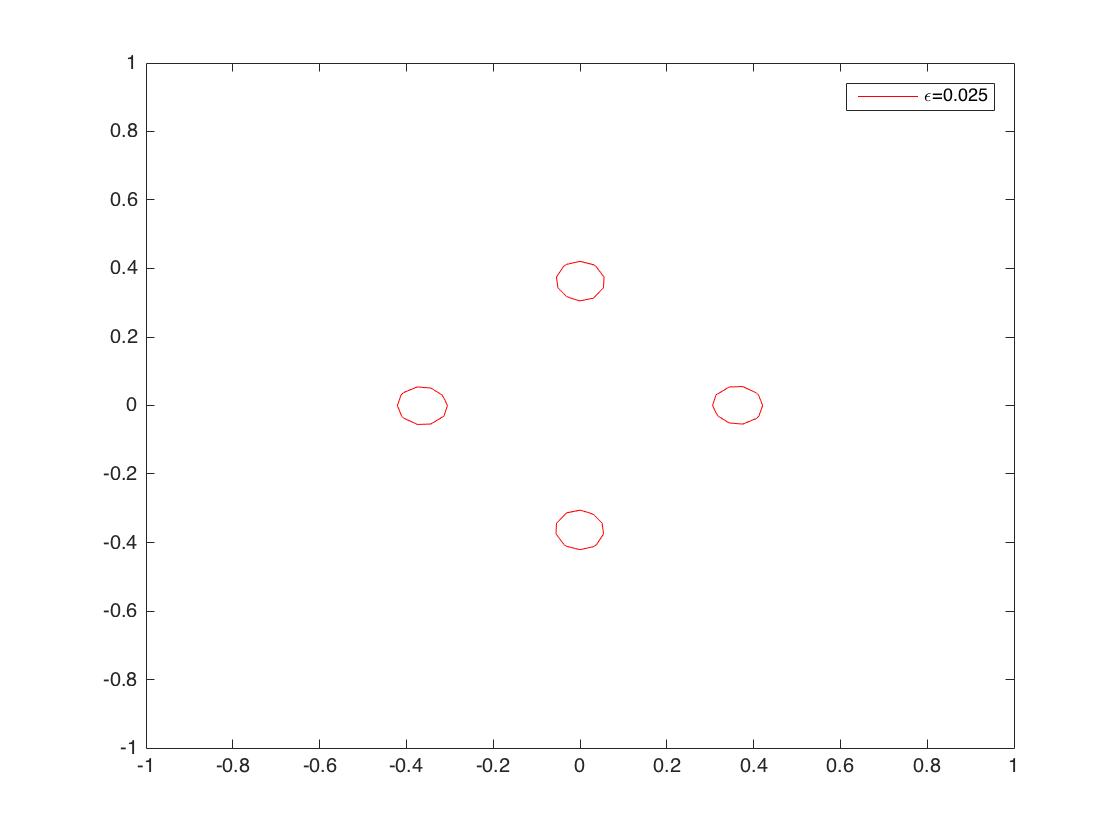}

\caption{Test 2 Snapshots of the zero-level set of $u^{\epsilon,h,k}$ at time $t=0,
5.5\times 10^{-3}, 1.1\times 10^{-2}, 1.65\times 10^{-2}, 2.2\times 10^{-2}, 2.75\times 10^{-2}$  and $\epsilon=0.125$.}\label{figure}
\end{figure}

Figure 5.2 displays six snapshots of the zero-level set of the numerical solution $u^{\epsilon,h,k}$ with $\epsilon=0.025$.
Similarly, we observe that as $\epsilon$ is small enough the zero-level set converges to the mean curvature flow $\Gamma_t$
as time goes on.

\vspace{10pt} \textbf{ Acknowledgment:}\ {\small The authors would like to express sincere thanks to Dr. Yukun Li of the Ohio State University for introducing Allen-Cahn equation to the authors and for his many valuable discussions and suggestions.}


\end{document}